\documentclass[10pt,reqno]{amsart}
\usepackage{a4wide}
\usepackage[mathscr]{eucal}
\usepackage{amssymb}
\textwidth=125mm  \textheight=195mm \DeclareMathOperator{\C}{C}
\DeclareMathOperator{\M}{M} \DeclareMathOperator{\Mor}{Mor}
\DeclareMathOperator{\B}{B} \DeclareMathOperator{\Rep}{Rep}

\def\build#1_#2^#3{\mathrel{\mathop{\kern 0pt#1}\limits_{#2}^{#3}}}

\newtheorem{twr}{Theorem}[section]
\newtheorem{lem}[twr]{Lemma}
\newtheorem{stwr}[twr]{Proposition}

\theoremstyle{definition}
\newtheorem{defin}[twr]{Definition}

\newtheorem{uw}[twr]{Remark}
\newcommand{\ir}{{\rm i}}
\newcommand{\id}{{\rm id}}
\newcommand{\rf}[1]{(\ref{#1})}
\newcommand{\aff}{\mathbin\eta}
\newcommand{\Ph}{{\rm Phase}}
\newcommand{\af}{\hspace*{.1cm}\eta \hspace*{.05cm}}
\begin{document}
\title[]{Rieffel deformation via crossed products}
\author{P.~Kasprzak}
\address{Department of Mathematical Methods in Physics\\
Faculty of Physics\\
Warsaw University} \email{pawel.kasprzak@fuw.edu.pl}
\thanks{The research was  supported by the KBN under grant 115/E-343/SPB/6.PR UE/DIE 50/2005 – 2008.}
\begin{abstract} We start from Rieffel data $(A,\Psi,\rho)$, where $A$ is a
$\C^*$-algebra, $\rho$ is an action of an abelian group $\Gamma$
on $A$ and $\Psi$ is a $2$-cocycle on the dual group. Using
Landstad theory of crossed product we get a deformed
$\C^*$-algebra $A^\Psi$. In the  case of $\Gamma=\mathbb{R}^n$ we
obtain a very simple proof of invariance of $\mathcal{K}$-groups
under the deformation. In the general case we also get a very
simple proof that nuclearity is preserved under the deformation.
We show how our approach leads
 to quantum groups and investigate their duality.
The general theory is illustrated by an example of the deformation
of $SL(2,\mathbb{C})$. A description of it, in terms of
noncommutative coordinates
$\hat\alpha,\hat\beta,\hat\gamma,\hat\delta$, is given.
\end{abstract}

\maketitle \tableofcontents\newpage

\begin{section}{Introduction.} In \cite{Rf1} Rieffel described the
method of deforming of $\C^*$-algebras known today as the
\emph{Rieffel deformation}. Having an action of $\mathbb{R}^d$ on
a $\C^*$-algebra $A$ and a skew symmetric operator $J:
\mathbb{R}^d\mapsto\mathbb{R}^d$ Rieffel defined a new product
that gave rise to the deformed $\C^*$-algebra $A^J$. In \cite{Rf2}
the Rieffel deformation was applied to the $\C^*$-algebra of
continuous functions vanishing at infinity
 on a Lie group $G$.
An action of $\mathbb{R}^n$ was constructed using the left and
right shifts along a fixed abelian Lie subgroup $\Gamma$.  Having
the deformed $\C^*$-algebra Rieffel introduced a comultiplication,
 a coinverse and a counit, showing that it is a locally compact quantum
group.

 M. Enock and
 L.  Vainerman in \cite{V} gave a method of deforming of the dual object
 associated with the locally compact group $G$ that is $(\C_r^*(G),\Hat\Delta)$
where $\C_r^*(G)$ is the reduced group $\C^*$-algebra and
$\Hat\Delta$ is the canonical comultiplication on it.
 Using an  abelian subgroup $\Gamma\subset G$  and a  $2$-cocycle $\Psi$
on the Pontryiagin  dual group $\Hat\Gamma$ they twisted the
canonical comultiplication $\Hat\Delta$ on the reduced group
$\C^*$-algebra $\C^*_r(G)$ obtaining a new quantum group. They
also presented  a formula for a multiplicative unitary and
described a Haar measure for this new quantum group.

The existence of these two methods of deforming of objects related
to a group $G$ prompts the question about the relations between
them. In this paper it is shown that they are dual versions of the
same mathematical procedure. Let us note that the deformation
framework of Enock and Vainerman is in a sense more general than
the one of Rieffel: instead of a skew symmetric matrix on
$\mathbb{R}^n$ they use a $2$-cocycle $\Psi$ on $\Hat\Gamma$. This suggests that it should be possible to
perform the Rieffel deformation of a $\C^*$-algebra $A$ acted on
by an abelian group $\Gamma$ with a $2$-cocycle $\Psi$ on
$\Hat{\Gamma}$. A formulation of Rieffel deformation in that
context is one of the  results of this paper.

Let us briefly describe the contents of the whole paper. In the
next section we  revise the Landstad theory of crossed products.
We prove a couple of useful results that we could not find in the
literature. In Section \ref{Riefdef} we use the Landstad's theory
to give a new approach to the Rieffel deformation of
$\C^*$-algebras. In Section \ref{qg} we apply the Rieffel
deformation to locally compact groups. We show that
Enock-Vainerman's and Rieffel's approach give mutually dual,
locally compact quantum groups. Moreover, a formula for a Haar
measure on a quantized algebra of functions is given.
 In the last section we use our scheme to deform  $SL(2,\mathbb{C})$.
The subgroup $\Gamma$ consists of diagonal matrices. We show that
the deformed $\C^*$-algebra $A$ is generated in the sense of
Woronowicz by four affiliated elements
$\hat{\alpha},\hat{\beta},\hat{\gamma},\hat{\delta}$ and give a
detailed description of the commutation relations they satisfy.
Moreover, we show that the comultiplication
$\Delta^\Psi\in\Mor(A;A\otimes A)$ acts on the generators in the
standard way:
\begin{equation}
 \label{introduction}
\begin{array}{cc}
  \Delta^\Psi(\hat{\alpha})= \hat{\alpha}\otimes \hat{\alpha}\dotplus \hat{\beta}\otimes \hat{\gamma}&
  \Delta^\Psi(\hat{\beta})= \hat{\alpha}\otimes \hat{\beta}\dotplus \hat{\beta}\otimes \hat{\delta} \\
   \Delta^\Psi(\hat{\gamma})= \hat{\gamma}\otimes \hat{\alpha}\dotplus \hat{\delta}\otimes \hat{\gamma}&
   \Delta^\Psi(\hat{\delta})= \hat{\delta}\otimes \hat{\delta}\dotplus \hat{\gamma}\otimes \hat{\beta}.\end{array}
\end{equation}
Throughout the paper we will freely use the language of
$\C^*$-algebras and the theory of locally compact quantum groups.
For the notion of multipliers, affiliated elements, algebras
generated by a family of affiliated elements and morphism of
$\C^*$-algebras we refer the reader to \cite{W4} and \cite{W10}.
For the theory of locally compact quantum groups we refer to
\cite{kv} and \cite{MNW}. 

Some remarks about the notation. For a subset $X$ of a Banach
space $B$, $X^{\rm cls}$ denotes the closed linear span of $X$.
Let $A$ be a $\C^*$-algebra and $A^*$ be its Banach dual. $A^*$ is
an $A$-bimodule where for $\omega\in A^*$ and $b,b'\in A$ we
define $b\cdot\omega\cdot b'$ by the formula:
\[b\cdot\omega\cdot b'(a)=\omega(b'ab)\]
for any $a\in A$.
\end{section}
\begin{section}{Landstad theory of crossed products.}\label{rozdzilk}
 Let us start this section with a definition of
$\Gamma$-product. For a  detailed treatment of this notion see
\cite{Pe}.
 \begin{defin}Let $\Gamma$ be a locally compact abelian group,
 $\Hat\Gamma$ its Pontryagin dual,
  $B$  a $\C^*$-algebra,
  $\lambda$  a homomorphism of
$\Gamma$ into the unitary group of $\,\M(B)$ continuous in the
strict topology of $\,\M(B)$ and let $\hat{\rho}$ be a continuous
action
 of $\hat{\Gamma}$ on $B$. The triple
$(B,\lambda,\hat{\rho})$ is called a $\Gamma$-product if:
\[\hat{\rho}_{\hat{\gamma}}(\lambda_\gamma)=\langle\hat{\gamma},\gamma\rangle\lambda_\gamma\]
for any $\hat\gamma\in\Hat\Gamma$ and $\gamma\in\Gamma$.
\end{defin}
The unitary representation $\lambda:\Gamma\rightarrow \M(B)$ gives
rise to a morphism of $\C^*$-algebras $\lambda\in
\Mor\bigl(\C^*(\Gamma);B\bigr)$. Identifying $\C^*(\Gamma)$ with $
\C_\infty( \hat{\Gamma})$ via the Fourier transform, we get a
morphism $\lambda\in \Mor\bigl( \C_\infty( \hat{\Gamma});B\bigr)$.
Let $\tau_{\hat{\gamma}}\in {\rm Aut}\bigl( \C_\infty(
\hat{\Gamma})\bigr)$ denote the shift automorphism:\[
\tau_{\hat{\gamma}}(f)(\hat{\gamma}')=f(\hat{\gamma}'+\hat{\gamma})
\mbox{  for all  } f\in  \C_\infty( \hat{\Gamma}).\] It is easy to
see that $\lambda$ intertwines the action $\hat{\rho}$ with
$\tau$:
\begin{equation}\label{int1}
\lambda(\tau_{\hat{\gamma}}(f))=\hat{\rho}_{\hat{\gamma}}(\lambda(f))
\end{equation} for any $f\in \C_\infty( \hat{\Gamma})$. The following lemma seems to be known but we could
not find any reference.
\begin{lem} \label{iota}  Let $(B,\lambda,\hat{\rho})$ be a $\Gamma$-product.
Then the morphism $\lambda\in \Mor( \C_\infty(
\hat{\Gamma});B)$ is injective.
\end{lem}
\begin{proof}  The kernel of the morphism $\lambda$ is an ideal in
$\C_\infty( \hat{\Gamma})$ hence it is contained in a maximal ideal.
Therefore there exists $\hat{\gamma}_0$ such that
$f(\hat{\gamma}_0)=0$ for all $f\in \ker\lambda$. Equation
\eqref{int1} implies that $\ker\lambda$ is $\tau$ invariant. Hence
$f(\hat{\gamma}_0+\hat{\gamma})=0$ for all $\hat{\gamma}$. This
shows that $f=0$ and $\ker\lambda=\{0\}$.\end{proof}

In what follows we usually treat a $\C^*$- algebra
$\C_\infty(\Hat\Gamma)$ as a subalgebra of $\M(B)$ and we will not
use the embedding $\lambda$  explicitly.
\begin{defin}\label{invL}
Let $(B,\lambda,\hat{\rho})$ be a $\Gamma$-product and  $x\in
\M(B)$. We say that $x$
 satisfies the Landstad  conditions if:
\begin{equation}\label{lc1}\left\{
\begin{array}{l}
 (i)\,\,\,\,\,\,\,\hat{\rho}_{\hat{\gamma}}(x)=x \mbox{ for all }\hat{\gamma}\in
\hat{\Gamma};\\(ii)\,\,\mbox{ the map }\Gamma\ni \gamma\mapsto
\lambda_\gamma x \lambda_\gamma^*\in\M(B) \mbox{ is norm continuous;
}\\(iii)\,\, fxg\in B \mbox{ for all } f,g\in
  \C_\infty( \hat{\Gamma}).
\end{array}\right.
\end{equation}
\end{defin}
 In computations it is useful to smear unitary
 elements
 $\lambda_\gamma\in \M(B)$  with a function $h\in
 L^1(\Gamma)$:
 \[\lambda_h=\int_\Gamma h(\gamma)\lambda_\gamma d\gamma \in
\M(B).\] Note that $\lambda_h\in \M(B)$ coincides with the Fourier
transform of $h$: $\mathfrak{F}(h)\in  \C_\infty( \hat{\Gamma})$.
In the original form of Definition \ref{invL} given by Landstad
the third condition had the form:
\begin{equation}\label{lc2} \lambda_fx,\,\, x\lambda_f\in B\,\,
{\rm for\,\, all} \,\, f,g\in
 L^1(\Gamma).\end{equation}   Our
 conditions are simpler to check, which turns out to be useful in
 the example considered at the end of the paper. As shown below
both definitions of invariants are in fact equivalent. The
argument  is very similar to the one given in \cite{PuszSoltan}
which shows that the third Landstad condition can be also replaced
by \[ \lambda_fx\in B\,\,\, {\rm for\,\, all} \,\, f\in
 L^1(\Gamma).\] Assume that
$x\in \M(B)$ satisfies the Landstad conditions \eqref{lc1}. Choose
$\varepsilon\geq 0$ and a function $f\in L^1(\Gamma)$. By the
second Landstad condition we can find a finite volume neighborhood
$ \mathcal{O}$ of the neutral element $e\in \Gamma$ such that:
\begin{equation}\label{fvol} \|\lambda_\gamma x-x\lambda_\gamma\|\leq\varepsilon
\,\,\,\,{\rm for}\,\,{\rm all}\,\,\,\, \gamma\in \mathcal{O}
.\end{equation} Let $\chi_{\mathcal{O}}$ denote the normalized
characteristic function of $\mathcal{O}\subset \Gamma$:
\begin{displaymath}
\chi_{\mathcal{O}}(\gamma)=\left\{\begin{array} {ll} \frac{1}{ vol
\mathcal{O}}&\,\,\, {\rm if}\,\,\, \gamma\in\mathcal{O}\\
0&\,\,\,{\rm if}\,\,\, \gamma\in \Gamma\setminus \mathcal{O}.
\end{array}\right.
\end{displaymath}
Then by \eqref{fvol} we have: \begin{equation}\label{fvol1}
\|\lambda_{\chi_{\mathcal{O}} }x-x \lambda_{\chi_{\mathcal{O}} }\|
 \leq\varepsilon.\end{equation} If necessary, we can choose  a smaller neighborhood
 and assume also that:
 \begin{equation}\label{fvol2}\|\lambda_f\lambda_{\chi_{\mathcal{O}} }
-\lambda_f\|\leq\varepsilon.\end{equation}
 The calculation below is self-explanatory
 \[\lambda_fx=\lambda_fx-\lambda_f\lambda_{\chi_{\mathcal{O}} } x+\lambda_f\lambda_{\chi_{\mathcal{O}} } x=
 (\lambda_fx-\lambda_f\lambda_{\chi_{\mathcal{O}} } x)
 +(\lambda_f\lambda_{\chi_{\mathcal{O}} } x-\lambda_fx\lambda_{\chi_{\mathcal{O}} } )
 +\lambda_fx\lambda_{\chi_{\mathcal{O}} }\]
 and together with \eqref{fvol1} and \eqref{fvol2} shows that:
  \[\|\lambda_fx-\lambda_fx\lambda_{\chi_{\mathcal{O}} }\|\leq \varepsilon(\|x\|+\|\lambda_f\|).\]
  Hence we can
 approximate $\lambda_fx$ by elements of the form
 $\lambda_fx\lambda_{\chi_{\mathcal{O}} }\in B$. This shows that
$\lambda_fx\in B$. A similar argument proves the second inclusion:
$x\lambda_f\in B$.

The set of elements satisfying Landstad's conditions is a
$\C^*$-algebra. We shall call it the Landstad algebra and denote
it by $A$. It follows from Definition \ref{invL}, that if $a\in A$
then $\lambda_\gamma a\lambda_\gamma^*\in A$  and the map
$\Gamma\in \gamma\mapsto \lambda_\gamma a\lambda_\gamma^*\in A$ is
norm continuous. An action of $\Gamma$ on $A$ defined in this way
will be denoted by $\rho$ .

 It can be shown that the embedding of $A$ into $\M(B)$ is a
morphism of $\C^*$-algebras (c.f. \cite{LOP}, Section 2). Hence
the multipliers algebra $\M(A)$ can also be  embedded into
$\M(B)$. Let $x\in \M(B)$. Then $x\in \M(A)$ if and only if it
satisfies the following two conditions:
\begin{equation}\label{lcm}\left\{
\begin{array}{l}
(i)\,\,\,\,\hat{\rho}_{\hat{\gamma}}(x)=x \mbox{ for all
}\hat{\gamma}\in \hat{\Gamma};\\(ii) \mbox{ for all } a\in A, \mbox{
the map } \\\,\,\,\,\,\,\,\,\,\,\,\Gamma\ni\gamma\mapsto
\lambda_\gamma x \lambda_\gamma^*a\in\M(B)\\\quad\,\,\, \mbox{ is
norm continuous.}
\end{array}\right.
\end{equation} Note that the first and the second condition of
\rf{lc1} imply conditions \rf{lcm}.

 Examples of $\Gamma$-products can be obtained via the
crossed-product construction. Let $A$ be a $\C^*$-algebra with an
action $\rho$ of $\Gamma$ on $A$. There exists the standard action
$ \hat{\rho}$ of the group $\hat{\Gamma}$ on $A\rtimes_\rho
\Gamma$ and a unitary representation
$\lambda_\gamma\in\M(A\rtimes_\rho \Gamma)$ such that the triple
$(A\rtimes_\rho \Gamma,\lambda_\gamma,\hat{\rho})$ is a
$\Gamma$-product. It turns out that all $\Gamma$-products
$(B,\lambda, \hat{\rho})$ are crossed-products of the Landstad
algebra $A$ by the action $\rho$ implemented by $\lambda$. The
following theorem is due to Landstad (Theorem 7.8.8, \cite{Pe}):
\begin{twr}\label{gspace} A triple  $(B,\lambda,\hat{\rho})$ is a
$\Gamma$-product if and only if there is a $\C^*$-dynamical system
$(A,\Gamma,\rho)$ such that $B=A\rtimes _\rho \Gamma$. This system
is unique up to isomorphism and $A$ consist of the elements in
$\M(B)$ that satisfy Landstad conditions while
$\rho_\gamma(a)=\lambda_\gamma a \lambda_\gamma^*$.
\end{twr}
\begin{uw}\label{uw1}\rm
The main problem in the proof of the above theorem is to show that
the Landstad algebra is not small. It is solved by integrating
the action $\hat{\rho}$ over the dual group. More precisely, we
say that an element $x\in \M(B)_+$ is $\hat{\rho}$ - integrable if
there exists $y\in \M(B)_+$ such that
\[\omega(y)=\int d\hat{\gamma}\,\omega(\hat{\rho}_{\hat{\gamma} }(x))\]
for any $\omega\in \M(B)^*_+$. We denote $y$ by $\mathfrak{E}(x)$.
If $x\in \M(B)$ is not positive then we say that it is
$\hat{\rho}$\,-integrable if it can be written as a linear
combination of positive $\hat{\rho}$\,-integrable elements. The
set of $\hat{\rho}$ - integrable elements will be denoted by
$D(\mathfrak{E})$. The averaging procedure induces  a  map:
\[\mathfrak{E}:D(\mathfrak{E})\mapsto \M(B).\]
It can be shown that for a large class of $x\in D(\mathfrak{E})$,
$\mathfrak{E}(x)$ is an element of the Landstad algebra $A$.  This
is the case for $f_1bf_2$ where $b\in B$ and $f_1,f_2\in  \C_\infty(
\hat{\Gamma})$ are square integrable. Moreover the map:
\begin{equation}\label{conav}B\ni b\mapsto \mathfrak{E}(f_1 bf_2)\in \M(B)\end{equation}
 is continuous with the following estimate for norms
\begin{equation}\label{conav1} \|\mathfrak{E}(f_1
bf_2)\|\leq\|f_1\|_2\|b\|\|f_2\|_2\end{equation} where $\|\cdot\|_2$
is the $L^2$-norm. Furthermore, we have
\begin{equation} \label{gsl}\left\{\mathfrak{E}(f_1bf_2) :\,\,b\in
B,\,\,f_1,f_2\in \C_\infty( \hat{\Gamma})\cap
L^2(\hat{\Gamma})\right \}^{{\rm cls}}=A.\end{equation} The last
equality  was not proven in \cite{Pe}. We shall need it at some
point so let us give a proof here. Let $a\in A$ and
$f_1,f_2,f_3,f_4$ be continuous, compactly supported functions on
$\Gamma$. Consider an element
$x=\lambda_{f_1}\lambda_{f_2}a\lambda_{f_3}\lambda_{f_4}$. Clearly
$x=\mathfrak{F}(f_1)\mathfrak{F}(f_2)a\mathfrak{F}(f_3)\mathfrak{F}(f_4)$
and  $\mathfrak{F}(f_1),\ldots,\mathfrak{F}(f_4)\in
L^2(\Hat\Gamma)$, hence by \eqref{conav1} $x\in D(\mathfrak{E})$.
We compute:
\begin{eqnarray*}
\mathfrak{E}(x)&\hspace*{-0,25cm}=&\hspace*{-0,25cm}\int
d\hat{\gamma}\,\hat{\rho}_{\hat{\gamma}
}(\lambda_{f_1}\lambda_{f_2}a\lambda_{f_3}\lambda_{f_4})\\
&\hspace*{-0,25cm}=&\hspace*{-0,25cm}\int
d\hat{\gamma}\,\hat{\rho}_{\hat{\gamma} }\left(\int
d\gamma_1d\gamma_2d\gamma_3d\gamma_4\,f_1(\gamma_1)f_2(\gamma_2)\lambda_{\gamma_1+\gamma_2}
a\lambda_{\gamma_3+\gamma_4}f_3(\gamma_3)f_4(\gamma_4)\right)\\
&\hspace*{-0,25cm}=&\hspace*{-0,25cm}\int d\hat{\gamma}\int
d\gamma_1d\gamma_2d\gamma_3d\gamma_4\Bigl(\langle\hat{\gamma},\gamma_1+\gamma_2+\gamma_3+\gamma_4\rangle
f_1(\gamma_1)f_2(\gamma_2)\\
&&\hspace*{4,5cm}\times\lambda_{\gamma_1+\gamma_2}
a\lambda_{\gamma_3+\gamma_4}f_3(\gamma_3)f_4(\gamma_4)\Bigr).
\end{eqnarray*}
Using properties of the Fourier transform we obtain :
\begin{equation}\label{apx} \mathfrak{E}(x)=\int
d\gamma_1d\gamma_2d\gamma_3\,\rho_{\gamma_1+\gamma_2}
(a)f_1(\gamma_1)f_2(\gamma_2)f_3(\gamma_3)f_4(-\gamma_1-\gamma_2-\gamma_3).
\end{equation} If  $f_1,f_2,f_3$ approximate the Dirac delta
function and $f_4(0)=1$ then using \eqref{apx} we see that elements
$\mathfrak{E}(\lambda_{f_1}\lambda_{f_2}a\lambda_{f_3}\lambda_{f_4})$
approximate $a$ in norm. This  proves \rf{gsl}.
\end{uw}
 The following lemma  is simple but very useful:
\begin{lem}\label{dsw} Let $ (B,\lambda,\hat{\rho})$ be a
$\Gamma$-product, A its Landstad algebra, $\rho$ an action of
$\Gamma$ on $A$ implemented by $\lambda$ and $ \mathcal{V}\subset
A$  a subset of the Landstad algebra which is invariant under the
action $\rho$ and such that $\left(\C_{\infty}(
\hat{\Gamma})\mathcal{V} \C_{\infty}( \hat{\Gamma})\right)^{\rm
cls}=B$. Then $\mathcal{V}^{\,\rm cls}=A$.
\end{lem}
\begin{proof}This proof is similar to the proof of formula
\rf{gsl}. Let $f_1,f_2,f_3,f_4$ be  continuous, compactly
supported functions on $\Gamma$. Using \eqref{conav1}  and
\eqref{gsl} we get:
\begin{equation}\label{denlan} A=\{\mathfrak{E}(\lambda_{f_1}(\lambda_{f_2}v\lambda_{f_3})\lambda_{f_4}):
v\in \mathcal{V}\}^{\rm cls}.\end{equation} A simple calculation
shows that:
\[\mathfrak{E}(\lambda_{f_1}(\lambda_{f_2}v\lambda_{f_3})\lambda_{f_4})=\int
d\gamma_1d\gamma_2d\gamma_3\,\rho_{\gamma_1+\gamma_2}(v)f_1(\gamma_1)
f_2(\gamma_2)f_3(\gamma_3)f_4(-\gamma_1-\gamma_2-\gamma_3).\] Note
that the integrand \[\Gamma\times \Gamma\ni
(\gamma_1,\gamma_2,\gamma_3)\mapsto
\rho_{\gamma_1+\gamma_2}(v)f_1(\gamma_1)f_2(\gamma_2)f_3(\gamma_3)f_4(-\gamma_1-\gamma_2-\gamma_3)\in
\mathcal{V}\] is a norm continuous, compactly supported function,
hence
\[ \mathfrak{E}(\lambda_{f_1}v\lambda_{f_2}\lambda_{f_3})\in
\mathcal{V}^{\rm cls} .\] From this and \eqref{denlan} we get
$\mathcal{V}^{\rm cls}=A$.
\end{proof}
The next proposition shows that morphisms of $\Gamma$-products
induce morphisms of their Landstad algebras. The below result can
be, to some extent, deduced from the results of paper \cite{KQR}.
\begin{stwr}\label{gmor1} Let $(B,\lambda,\hat{\rho})$ and $(B',\lambda',\hat{\rho}')$
be $\Gamma$-products and let $A$, $A'$ be Landstad algebras for
$B$ and $B'$ respectively. Assume that $\pi\in \Mor(B,B')$
satisfies:
\begin{itemize}
\item $\pi(\lambda_\gamma)=\lambda'_\gamma$;
\item $\pi\bigl( \hat{\rho}_{ \hat{\gamma}}(b)\bigr)=\hat{\rho}'_{
\hat{\gamma}}\bigl(\pi(b)\bigr)$.
\end{itemize}
Then $\pi(A)\subset \M(A')$ and $\pi|_{A}\in \Mor(A,A')$.
Moreover,  if $\pi(B)\subset B'$ then $\pi(A)\subset A'$. If
$\pi(B)= B$ then $\pi(A)= A'$.
\end{stwr}
\begin{proof}
We start by showing that $\pi(A)\subset \M(A')$. Let $a\in A$. Then
 \[\hat{\rho}'_{ \hat{\gamma}}(\pi(a))=\pi(\hat{\rho}_{
\hat{\gamma}}(a))=\pi(a).\] Hence $\pi(a)$ is $\hat{\rho}'$
invariant. Moreover the map: \[\Gamma\ni\gamma\mapsto
{\lambda'}_\gamma \pi(a){\lambda'}_\gamma^*=\pi(\lambda_\gamma a
\lambda_\gamma^*)\in \M(A')\] is norm continuous. This shows that
$\pi(a)$ satisfies the first and the second Landstad condition of
\eqref{lc1} which guaranties that $\pi(a)\in \M(A')$.

To prove that the homomorphism $\pi$ restricted to $A$ is in fact
a morphism from $A$ to $A'$ we have to check that the set
$\pi(A)A'$ is linearly dense in $A'$. We know that $\pi(B)B'$ is
linearly dense in $B'$. Using this fact in the last equality
below, we get
\begin{eqnarray*}\left(\C^*(\Gamma)\pi(A)A'\C^*(\Gamma)\right)^{\rm cls}
&\hspace*{-0,25cm}=&\hspace*{-0,25cm}\left(\pi(\C^*(\Gamma)A)A'\C^*(\Gamma)\right)^{\rm cls}\\
&\hspace*{-0,25cm}=&\hspace*{-0,25cm}\left(\pi(B)B'\right)^{\rm
cls}=B'.\end{eqnarray*} Moreover
\[\lambda'_\gamma \pi(a)a'{\lambda'}_\gamma^*=\pi(\lambda_\gamma a
\lambda_\gamma^*)\lambda'_\gamma a'{\lambda'}_\gamma^*,\] hence
the set $\pi(A)A'$ is $\rho'$-invariant (remember that $\rho'$ is
the action of $\Gamma$ implemented by $\lambda'$). This shows that
$\pi(A)A'$ satisfies the assumptions of Lemma \ref{dsw} and gives
the density of $\pi(A)A'$ in $A'$.

Assume now that $\pi(B)\subset B'$. Let $a\in A$ satisfy
 Landstad conditions \rf{lc1}.
Then as was shown at the beginning of the proof,  $\pi(a)$
satisfies the first and the second Landstad condition. Moreover
$f\pi(a)g=\pi(fag)\in B' \mbox{  for all  } f,g\in \C_{\infty}(
\hat{\Gamma})$, hence $\pi(a)$  also satisfies the third Landstad
condition. Therefore $\pi(a)\in A'$.

If $\pi(B)=B'$ then the equality
\begin{equation}\label{avr}\mathfrak{E}(
f_1\pi(b)f_2)=\pi(\mathfrak{E}( f_1bf_2))\end{equation} and property
\rf{gsl} shows that $\pi(A)=A'$.

To prove \rf{avr} take $\omega\in \M(B')^*$. Then
\begin{eqnarray*} \omega(\pi(\mathfrak{E}(
f_1bf_2)))&\hspace*{-0,25cm}=&\hspace*{-0,25cm}\omega\circ\pi(\mathfrak{E}(
f_1bf_2)))=\int d \hat{\gamma}\,\omega\circ\pi (
\hat{\rho}_{\hat{\gamma}}(f_1bf_2))\\&\hspace*{-0,25cm}=&\hspace*{-0,25cm}\int
d \hat{\gamma}\,\omega(
\hat{\rho}_{\hat{\gamma}}(f_1\pi(b)f_2))=\omega\left(\int d
\hat{\gamma}\,
\hat{\rho}_{\hat{\gamma}}(f_1\pi(b)f_2)\right)\\&\hspace*{-0,25cm}=&\hspace*{-0,25cm}\omega(\mathfrak{E}(
f_1\pi(b)f_2)).
\end{eqnarray*}
Hence $\pi(\mathfrak{E}( f_1bf_2))=\mathfrak{E}( f_1\pi(b)f_2)$.
\end{proof}
Let $\Gamma'$ be an abelian locally compact group and
$\phi:\Gamma\mapsto \Gamma'$  a continuous homomorphism. For
$\hat{\gamma}'\in \hat{\Gamma}'$ we set
$\phi^T(\hat{\gamma}')=\hat{\gamma}'\circ \phi\in \hat{\Gamma}$.
The map
\[ \phi^T: \hat{\Gamma}'\mapsto \hat{\Gamma},\quad
\phi^T(\hat{\gamma}')=\hat{\gamma}'\circ \phi\] is a continuous
group homomorphism called the dual homomorphism. We have a version
of Proposition \ref{gmor1} with two different groups.
\begin{stwr}\label{gmor3} Let $(B,\lambda,\hat{\rho})$ be a $\Gamma$-product,
$(B',\lambda',\hat{\rho}')$ a $\Gamma'$-product,
$\phi:\Gamma\mapsto\Gamma'$ a surjective continuous homomorphism
and $\phi^T:\hat{\Gamma}'\mapsto\hat{\Gamma}$ the dual
homomorphism. Assume that $\pi\in \Mor(B,B')$ satisfies:
\begin{itemize}
\item $\pi(\lambda_{\gamma})=\lambda'_{\phi(\gamma)}$
\item $\hat{\rho}'_{\hat{\gamma}'}\bigl(\pi(b)\bigr)=\pi\bigl(\hat{\rho}_{
\phi^T(\hat{\gamma}')}(b)\bigr).$ \end{itemize}  Then
$\pi(A)\subset \M(A')$ and $\pi|_{A}\in {\rm Mor(A,A')}$.
Moreover, if $\pi(B)\subset B'$ then $\pi(A)\subset A'$. If
$\pi(B)= B$ then $\pi(A)= A'$.
\end{stwr}

Let $\pi\in \Mor(B,B')$ be a morphism of $\C^*$-algebras
satisfying the assumptions of Proposition \ref{gmor1} such that
$\pi(B)=B'$. We have an exact sequence of
$\C^*$-algebras:\begin{equation}\label{esi}0\rightarrow
\ker\pi\rightarrow B\build{\rightarrow}_{}^{\pi }B'\rightarrow
 0.\end{equation}  The $\C^*$-algebra $\ker\pi$ has a canonical $\Gamma$-product structure. Indeed,
consider a morphism $\alpha\in \Mor(B;\ker\pi)$ associated with
the ideal $\ker\pi\subset B$:\[\alpha(b)j=bj\] where $b\in\B$ and
$j\in\ker\pi$. Note that
\begin{equation}\label{alker}
\alpha(b)=b \mbox{ for any } b\in\ker\pi\subset B.
\end{equation}
 For all $\gamma\in\Gamma$ we set
$\tilde{\lambda}_\gamma= \alpha(\lambda_\gamma)\in \M(\ker\pi)$. The
map $\Gamma\ni\gamma\mapsto\tilde{\lambda}_\gamma\in \M(\ker\pi)$
  is a strictly continuous representation of
$\Gamma$ on $\ker\pi$.
 Moreover   $\ker \pi$
 is invariant under the action $\hat{\rho}$. The restriction of $\hat\rho$ to $\ker\pi$
will also be  denoted by  $\hat\rho$.
  It is easy to check that
  $\hat\rho_{\hat{\gamma}}\bigl(\tilde{\lambda}_\gamma\bigr)=\langle \hat{\gamma},\gamma\rangle
 \tilde{\lambda}_\gamma$
 which shows that the triple
 $\bigl(\ker\pi,\tilde{\lambda},\hat{\rho}\bigr)$ is a $\Gamma$-product.
Let $\mathcal{I}$, $A$, $A'$ be Landstad algebras for the
$\Gamma$-products $\bigl(\ker
\pi,\tilde{\lambda},\hat{\rho}\bigr)$, $(B,\lambda,\hat{\rho})$,
$(B',\lambda',\hat{\rho}')$ respectively. Our objective is to show
that the exact sequence \rf{esi} induces an exact sequence of
Landstad algebras:
\begin{equation}\label{esi1}0\rightarrow \mathcal{I}\rightarrow A\rightarrow
A'\rightarrow
 0.\end{equation}  Let $\bar{\pi}\in \Mor( A; A')$ denote a
morphism of Landstad algebras induced by $\pi$. We assumed that
$\pi$ is surjective, hence by Proposition \ref{gmor1}
$\bar\pi(A)=A'$ and we have an exact sequence of $\C^*$-algebras:
\[0\rightarrow \ker \bar{\pi}\rightarrow
A\build{\rightarrow}_{}^{\bar{\pi} }A'\rightarrow
 0.\]
It is easy to check that the morphism $\alpha\in \Mor(B;\ker \pi)$
satisfies the assumptions of Proposition \ref{gmor1},
 hence $\alpha(A)\subset \M(\mathcal{I})$.
 If we show that $\alpha$ restricted to $ \ker \bar{\pi}$ identifies
it with the Landstad algebra $\mathcal{I}$, then the existence of
the exact sequence \rf{esi1} will be proven. There are two
conditions to be checked :
 \begin{itemize}
 \item[(i)]$\alpha( \ker \bar{\pi})=\mathcal{I}$;
 \item[(ii)] if $x\in \ker \bar{\pi}$ and $\alpha(x)=0$  then
 $x=0$.
 \end{itemize}
Ad(i)
 Let $a\in\ker \bar{\pi}$ and $f\in
 \C_\infty(\hat{\Gamma})$.
 Then $af\in B\cap \ker \pi$,
 hence \begin{equation}\label{al3} \alpha(a)f=\alpha(af)=af\in\ker \pi,\end{equation}
 where  we used \eqref{alker}.
This shows that $\alpha(a)$ satisfies the third Landstad condition
for $\Gamma$-product $\bigl(\ker
\pi,\tilde{\lambda},\hat{\rho}\bigr)$.
 As in Proposition \ref{gmor1} we check that $\alpha(a)$ also satisfies
 the first and the second Landstad condition,
 hence $\alpha(\ker \bar{\pi})\subset \mathcal{I}$.
Furthermore,  $\mathfrak{E}(f_1bf_2)\in {\ker \bar{\pi}}$ for all
$b\in \ker\pi$, and we have \[ \alpha(\mathfrak{E}(f_1bf_2)=
\mathfrak{E}(f_1\alpha(b)f_2))=\mathfrak{E}(f_1bf_2).\]
 Using \eqref{gsl} we see that $\alpha( \ker
 \bar{\pi})=\mathcal{I}$.\\
Ad(ii) Assume that $a\in\ker \bar{\pi}$ and $\alpha(a)=0$. Note
that $af\in\ker \pi$ for any $f\in \C_\infty(\hat{\Gamma})$. Using
\eqref{alker} we get  $af=\alpha(af)=\alpha(a)\alpha(f)=0$. Hence
$af=0$ for any $f\in\C_\infty(\Hat\Gamma)$, which implies that
$a=0$. We can summarize the above considerations in the following:
\begin{stwr} \label{exct}Let $(B,\lambda, \hat{\rho})$, $(B',\lambda',
\hat{\rho}')$ be $\Gamma$-products with Landstad algebras $A$,
$A'$ respectively, $\pi\in \Mor(B,B')$ a surjective morphism
intertwining $\hat{\rho}$ and $\hat{\rho}'$ such that
$\pi(\lambda_\gamma)=\lambda'_\gamma$. Let $(\ker
\pi,\tilde{\lambda},\hat{\rho})$ be the $\Gamma$-product described
after Proposition \ref{gmor3} and let $\mathcal{I}\subset\M({\ker
\pi})$ be its Landstad algebra. Then $\mathcal{I}$  can be
embedded into $A$ and we have a $\Gamma$-equivariant exact
sequence:
\[0\rightarrow\mathcal{ I}\rightarrow A\build{\rightarrow}_{}^{\bar{\pi}
}A'\rightarrow
 0\] where $\bar\pi=\pi|_A$.
 \end{stwr}
\end{section}
\begin{section}{Rieffel deformation of
$\C^*$-algebras.}\label{Riefdef}
\begin{subsection}{Deformation procedure.}
 Let $(B,\lambda, \hat{\rho})$ be a $\Gamma$-product.
A $2$-cocycle on the group $\hat\Gamma$ is a  continuous function
  $\Psi:\hat{\Gamma}\times
\hat{\Gamma}\rightarrow
 \mathbb{T}^1$ satisfying:
\begin{itemize}
\item[(i)] $\Psi(e,\hat\gamma)=\Psi(\hat\gamma,e)=1$ for all $\hat\gamma\in\Hat\Gamma$;
\item[(ii)]$\Psi(\hat{\gamma}_1,\hat{\gamma}_2+\hat{\gamma}_3)
\Psi(\hat{\gamma}_2,\hat{\gamma}_3)
=\Psi(\hat{\gamma}_1+\hat{\gamma}_2,\hat{\gamma}_3)\Psi(
\hat{\gamma}_1,\hat{\gamma}_2) $ for all
$\hat{\gamma}_1,\hat{\gamma}_2,\hat{\gamma}_3\in\Hat\Gamma$.
\end{itemize}
(For the theory of $2$-cocycles we refer to \cite{kl}.)

For $\hat{\gamma},\hat{\gamma}_1$ we set  $\Psi_{\hat{\gamma}}(
\hat{\gamma}_1)=
 \Psi( \hat{\gamma}_1,\hat{\gamma})$. It defines a family of functions  $\Psi_{\hat{\gamma}}:
 \hat{\Gamma}\mapsto
 \mathbb{T}^1$.
 Using the embedding $\lambda\in \Mor( \C_\infty(\hat{\Gamma});B)$
 we get a strictly continuous family of unitary elements
\begin{equation}\label{psig} U_{\hat{\gamma}}=\lambda(\Psi_{\hat{\gamma}})\in \M(B).\end{equation}
 The $2$-cocycle condition for $\Psi$ gives:
\begin{equation}\label{cocycl} U_{\hat{\gamma}_1+\hat{\gamma}_2}=\overline{\Psi(\hat{\gamma}_1,\hat{\gamma}_2)}U_{
\hat{\gamma}_1}\hat{\rho}_{{\hat\gamma}_1}(U_{\hat{\gamma}_2}).\end{equation}
\begin{twr}\label{defgam} Let $(B,\lambda,\hat{\rho})$ be a $\Gamma$-product and let $\Psi$
be a $2$-cocycle on $ \hat{\Gamma}$. For any $ \hat{\gamma}\in \hat{\Gamma}$ the map
\[\hat{\rho}^\Psi_{ \hat{\gamma}}:B\ni b\mapsto \hat{\rho}^\Psi_{ \hat{\gamma}}(b)=
U^*_{\hat{\gamma}}\hat{\rho}_{\hat{\gamma}}(b)U_{\hat{\gamma}}\in
B\]  is an automorphism of $\C^*$-algebra $B$. Moreover,
\[\hat{\rho}^\Psi:\hat{\Gamma} \ni \hat{\gamma}\mapsto \hat{\rho}^
\Psi_{\hat{\gamma}}\in {\rm Aut(B)}\] is a strongly continuous
action of $\,\,\hat{\Gamma}$ on $B$ and the triple
$(B,\lambda_\gamma, \hat{\rho}^\Psi)$ is a $\Gamma$-product.
\end{twr} \begin{proof} Using equation \rf{cocycl} we get
\begin{eqnarray*}
\hat{\rho}^\Psi_{\hat{\gamma}_1+\hat{\gamma}_2}(b)&\hspace*{-0,25cm}=&\hspace*{-0,25cm}
U^*_{\hat{\gamma}_1+\hat{\gamma}_2}\hat{\rho}_{\hat{\gamma}_1+\hat{\gamma}_2}(b)
U_{\hat{\gamma}_1+\hat{\gamma}_2}\\
&\hspace*{-0,25cm}=&\hspace*{-0,25cm}
\overline{\Psi(\hat{\gamma}_1,\hat{\gamma}_2)}U_{
\hat{\gamma}_1}^*\hat{\rho}_{{\gamma}_1}(U_{\hat{\gamma}_2})^*
\hat{\rho}_{\hat{\gamma}_1}(\hat{\rho}_{\hat{\gamma}_2}(b))\hat{\rho}_{{\gamma}_1}
(U_{\hat{\gamma}_2})
U_{\hat{\gamma}_1}\Psi(\hat{\gamma}_1,\hat{\gamma}_2)\\
&\hspace*{-0,25cm}=&\hspace*{-0,25cm}U_{
\hat{\gamma}_1}^*\hat{\rho}_{{\gamma}_1}(U_{\hat{\gamma}_2}^*\hat{\rho}_{\hat{\gamma}_2}(b)U_{\hat{\gamma}_2})U_{
\hat{\gamma}_1} =\hat{\rho}^\Psi_{\hat{\gamma}_1}(\hat{\rho}^\Psi_
{\hat{\gamma}_2}(b)).\end{eqnarray*} This shows that
$\hat{\rho}^\Psi$ is an action of $\hat{\Gamma}$ on $B$.
 Applying  $\hat{\rho}^\Psi$ to
$\lambda_\gamma$ we get:
\[\hat{\rho}^\Psi_{\hat{\gamma}}(\lambda_\gamma)=U^*_{\hat{\gamma}}\hat{\rho}_{
\hat{\gamma}}(\lambda_\gamma)U_{\hat{\gamma}}=\langle\hat{\gamma},\gamma\rangle
U^*_{\hat{\gamma}}\lambda_\gamma
U_{\hat{\gamma}}=\langle\hat{\gamma},\gamma\rangle\lambda_\gamma.\]
The last equality follows from commutativity of $\Gamma$. Hence the
triple $(B,\lambda_\gamma, \hat{\rho}^\Psi)$ is a $\Gamma$-product.
\end{proof}

The above theorem leads to the following procedure of  deformation of $\C^*$-algebras.
The data needed to perform the deformation is a
triple $(A,\rho,\Psi)$ consisting of a $\C^*$-algebra $A$, an
action $\rho$ of a  locally compact abelian group $\Gamma$ and a
$2$-cocycle $\Psi$ on $\Hat{\Gamma}$.  Such a triple is called
deformation data. The resulting $\C^*$-algebra will be denoted
$A^\Psi$. The procedure is carried out in three steps :
\begin{enumerate}
\item Construct the crossed product $B=A\rtimes_\rho\Gamma$. Let
$(B,\lambda,\Hat{\rho})$ be the standard $\Gamma$-product structure
of the crossed product.
\item  Introduce a $\Gamma$-product $\bigl(B,\lambda,\Hat{\rho}^\Psi\bigr)$
 as described in Theorem \ref{defgam}.
\item Let $A^\Psi$ be the Landstad algebra of the $\Gamma$-product
$\bigl(B,\lambda,\Hat{\rho}^\Psi\bigr)$.
\end{enumerate}
Note that $A^\Psi$ still carries an action $\rho^\Psi$ of $\Gamma$
given by
\[
\rho^\Psi_\gamma(x)=\lambda_\gamma{x}\lambda_\gamma^*.
\]
In this case it is not the formula defining the action itself, but
its domain of definition that changes under deformation. The
triple $(A^\Psi,\Gamma,\rho^\Psi)$ will be called a \emph{twisted
dynamical system}. The procedure of deformation described above is
called the Rieffel deformation.  Using Theorem \ref{gspace} we
immediately get
\begin{stwr}\label{eqgam} Let $(A,\rho,\Psi)$  be deformation data
and $(A^\Psi,\Gamma,\rho^\Psi)$ be the twisted dynamical system
considered above. Then
\[A\rtimes_\rho\Gamma=A^\Psi\rtimes _{\rho^\Psi}\Gamma.\]
\end{stwr}
In what follows we investigate the dependence of the Rieffel
deformation on the choice of a $2$-cocycle. Let
$f:\hat{\Gamma}\mapsto\mathbb{T}^1$ be a continuous function such
that $f(e)=1$. For all
$\hat{\gamma}_1,\hat{\gamma}_2\in\hat{\Gamma}$ we set
\[\partial f(\hat{\gamma}_1,\hat{\gamma}_2)=\frac{f(\hat{\gamma}_1+\hat{\gamma}_2)}
{f(\hat{\gamma}_1)f(\hat{\gamma}_2)}.\] One can check that the map
\[\partial f:\hat{\Gamma}\times \hat{\Gamma}\ni(\hat{\gamma}_1,\hat{\gamma}_2)\mapsto
 \frac{f(\hat{\gamma}_1+\hat{\gamma}_2)}{f(\hat{\gamma}_1)f(\hat{\gamma}_2)}
\in \mathbb{T}\] is a $2$-cocycle. $2$-cocycles of this form are
considered to be trivial. We say that a pair of $2$-cocycles $
\Psi_1, \Psi_2$ is in the same cohomology class if they differ by
a trivial $2$ cocycle: $\Psi_2=\Psi_1\partial f$.
\begin{twr}\label{cohcoc} Let $(A,\rho,\Psi)$ be deformation data, giving rise to
a Landstad algebra $A^\Psi$.
 Then the isomorphism class of the Landstad algebra $A^\Psi$ depends only on the cohomology class of $\Psi$.
\end{twr}
Theorem \ref{cohcoc} easily follows from the next two lemmas.
\begin{lem}\label{triv} Let $(A,\rho,\Psi)$ be deformation data with a trivial $2$-cocycle
$\Psi=\partial f$. Then $A$ and $A^\Psi$ are isomorphic. More
precisely, treating $f$ as an element of $\C^*$-algebra
$\M(A\rtimes_\rho\Gamma)$ we have $A^\Psi=\left\{faf^*:a\in
A\right\}$.
\end{lem}
\begin{proof}
Fixing the second  variable in $\Psi$ we get a family $\Psi_{
\hat{\gamma}}$ of the form
\begin{equation}\label{trivu}\Psi_{
\hat{\gamma}}=\overline{f(\hat{\gamma})}f^*\tau_{\hat{\gamma}}(f).\end{equation}
where
$\tau_{\hat{\gamma}}(f)(\hat{\gamma}')=f(\hat{\gamma}+\hat{\gamma}')$.
Let $U_{\hat{\gamma}}\in \M(A\rtimes_\rho\Gamma)$ be the unitary
element given by $\Psi_{ \hat{\gamma}}$ (c.f. \eqref{psig}). The
function $f$ can be embedded into $\M(A\rtimes_\rho \Gamma)$ and
using \rf{trivu} we get:
\begin{equation} \label{forU}
U_{\hat{\gamma}}=\overline{f(\hat{\gamma})}f^*\hat{\rho}_{\hat{\gamma}}(f).\end{equation}
Assume that $a\in A$. Then \begin{eqnarray*} \hat{\rho}^{\Psi}_{
\hat{\gamma}}(faf^*)&\hspace*{-0,25cm}=&\hspace*{-0,25cm}U^*_{\hat{\gamma}}\hat{\rho}_{
\hat{\gamma}}(faf^*)U_{\hat{\gamma}}\\
&\hspace*{-0,25cm}=&\hspace*{-0,25cm}U^*_{\hat{\gamma}}\hat{\rho}_{
\hat{\gamma}}(f)a\hat{\rho}_{
\hat{\gamma}}(f)^*U_{\hat{\gamma}}.\end{eqnarray*} Using equation
\rf{forU} we see that
\begin{eqnarray*}\hat{\rho}^{\Psi}_{
\hat{\gamma}}(faf^*)
&\hspace*{-0,25cm}=&\hspace*{-0,25cm}f(\hat{\gamma})f\hat{\rho}_{\hat{\gamma}}(f)^*\hat{\rho}_{
\hat{\gamma}}(f)a\hat{\rho}_{ \hat{\gamma}}(f)^*
\overline{f(\hat{\gamma})}f^*\hat{\rho}_{\hat{\gamma}}(f)
=faf^*\end{eqnarray*} which means that the element $faf^*$
satisfies the first Landstad condition for the $\Gamma$-product
$(A\rtimes_\rho\Gamma,\lambda,\hat{\rho}^{\Psi})$. It is easy to
check that it also satisfies the second and  third Landstad
condition, hence $fAf^*\subset A^\Psi$. An analogous reasoning
proves the opposite inclusion $fAf^*\supset A^\Psi$.
\end{proof}

 Let $\Psi_1,\Psi_2$ be a pair of  $2$-cocycles on $\Hat\Gamma$. Their product
$\Psi_1\Psi_2$ is also a $2$-cocycle. Let $(A,\Gamma,\rho)$ be a
dynamical system. The deformation data $(A,\rho,\Psi_1)$ gives
rise to the twisted dynamical system
$(A^{\Psi_1},\Gamma,\rho^{\Psi_1})$. Furthermore, the triple
$(A^{\Psi_1},\rho^{\Psi_1},\Psi_2)$ is deformation data which
gives rise to the $\C^*$-algebra $(A^{\Psi_1})^{\Psi_2}$. At the
same time, using the deformation data $(A,\rho,\Psi_1\Psi_2)$ we
can introduce the $\C^*$-algebra $A^{\Psi_1\Psi_2}.$
\begin{lem}\label{skladanie} Let $(A,\Gamma,\rho)$ be a $\Gamma$-product and
let $\Psi_1,\Psi_2$ be $2$-cocycles on the group $ \hat{\Gamma}$.
Let $(A^{\Psi_1})^{\Psi_2}$ be a $\C^*$-algebra constructed from
the deformation data $(A^{\Psi_1},\rho^{\Psi_1},\Psi_2)$ and let
$A^{\Psi_1\Psi_2}$ be a $\C^*$-algebra constructed from the
deformation data $(A,\rho,\Psi_1\Psi_2)$. Then
\[ A^{\Psi_1\Psi_2}\simeq(A^{\Psi_1})^{\Psi_2}.\]
\end{lem}
\begin{proof} The algebras $A^{\Psi_1\Psi_2}$ and $(A^{\Psi_1})^{\Psi_2}$
can be embedded into $\M(A\rtimes_\rho\Gamma )$: they are Landstad
algebras of the $\Gamma$-products $(A\rtimes_\rho \Gamma, \lambda,
\hat{\rho}^{\Psi_1\Psi_2})$ and $(A\rtimes_\rho \Gamma,
\lambda,({\hat{\rho}^{\Psi_1}})^{\Psi_2})$ respectively. Note that
$U^{\Psi_1\Psi_2}_{ \hat{\gamma}}=U^{\Psi_1}_{
\hat{\gamma}}U^{\Psi_2}_{ \hat{\gamma}}$, hence
\begin{eqnarray*}
\hat{\rho}^{\Psi_1\Psi_2}_{ \hat{\gamma}}(b)
&\hspace*{-0,25cm}=&\hspace*{-0,25cm}{U^{\Psi_1\Psi_2}}^*_{
\hat{\gamma}}\hat{\rho}_{ \hat{\gamma}}(b)U^{\Psi_1\Psi_2}_{
\hat{\gamma}}\\
&\hspace*{-,25cm}=&\hspace*{-,25cm}{U^{\Psi_2}}^*_{
\hat{\gamma}}{U^{\Psi_1}}^*_{ \hat{\gamma}} \hat{\rho}_{
\hat{\gamma}}(b)U^{\Psi_1}_{ \hat{\gamma}}U^{\Psi_2}_{
\hat{\gamma}}=({\hat{\rho}^{\Psi_1}})^{\Psi_2}_{
\hat{\gamma}}(b).
\end{eqnarray*} This shows that
$\hat{\rho}^{\Psi_1\Psi_2}=({\hat{\rho}^{\Psi_1}})^{\Psi_2}$ and
implies that the
  $(A\rtimes_\rho
\Gamma, \lambda, \hat{\rho}^{\Psi_1\Psi_2})$ and  $(A\rtimes_\rho
\Gamma, \lambda,({\hat{\rho}^{\Psi_1}})^{\Psi_2})$ are in fact the
same $\Gamma$-products. Therefore their Landstad algebras
coincide.
\end{proof}
\end{subsection}
\begin{subsection}{Functorial properties of the Rieffel deformation.}
Let $(B,\lambda, \hat{\rho})$ be a $\Gamma$-product, $\Psi$  a
$2$-cocycle on the dual group $ \hat{\Gamma}$ and $H$ a Hilbert
space. Using Theorem \ref{defgam} we introduce the twisted
$\Gamma$-product $(B,\lambda, \hat{\rho}^\Psi)$. Let
$A,A^\Psi\subset \M(B)$ be Landstad algebras of $(B,\lambda,
\hat{\rho})$ and $(B,\lambda, \hat{\rho}^\Psi)$ respectively and
$\pi\in {\rm Rep}(B;H)$  a representation of the $\C^*$-algebra
$B$. The representation of $B$ extends to multipliers $\M(B)$ and
can be restricted to $A$ and $A^\Psi$.
\begin{twr} \label{faithfulrep} Let $(B,\lambda, \hat{\rho})$, $(B,\lambda, \hat{\rho}^\Psi)$
be $\Gamma$-products considered above, $A$, $A^\Psi$ their
Landstad algebras and $\pi$ a representation of $\C^*$-algebra $B$
on a Hilbert space $H$. Then $\pi$ is faithful on $A$ if and only
if it is faithful on $A^\Psi$.
\end{twr}
\begin{proof} Assume that $\pi$ is faithful on $A$ and  let $a\in A^\Psi$
be such  that $\pi(a)=0$.
Invariance of $a$ with respect to the action $ \hat{\rho}^\Psi$
implies that
$\hat{\rho}_{\hat{\gamma}}(a)=U_{\hat{\gamma}}aU_{\hat{\gamma}}^*$.
Hence
\begin{equation}\label{faitrep}\hat{\rho}_{\hat{\gamma}}(f)U_{\hat{\gamma}}a
U_{\hat{\gamma}}^*\hat{\rho}_{\hat{\gamma}}(gf^*\lambda_{-\gamma})
=\hat{\rho}_{\hat{\gamma}}(fagf^*\lambda_{-\gamma}).\end{equation}
for all $f,g\in  \C_\infty( \hat{\Gamma})$. The element $a\in
A^\Psi$ belongs to $\ker \pi$ therefore
\[\pi\bigl(\hat{\rho}_{\hat{\gamma}}(f)U_{\hat{\gamma}}aU_{\hat{\gamma}}^*
\hat{\rho}_{\hat{\gamma}}(gf^*\lambda_{-\gamma})\bigr)=0.\]
Combining it with equation \eqref{faitrep} we obtain
$\pi\bigl(\hat{\rho}_{\hat{\gamma}}(fagf^*\lambda_{-\gamma})\bigr)=0$.

Assume now that $f,g\in L^2(\hat{\Gamma})\cap \C_\infty(
\hat{\Gamma})$. Let $\mathfrak{E}$ denote the averaging map with
respect to  undeformed action $\hat\rho$. Then
$fagf^*\lambda_{-\gamma}\in D(\mathfrak{E})$ and
$\mathfrak{E}(fagf^*\lambda_{-\gamma})=0$. Indeed, let $\omega\in
\B(H)_*$. Then
\begin{eqnarray*}
\omega\Bigl(\pi\bigl(\mathfrak{E}(fagf^*\lambda_{-\gamma})\bigr)\Bigr)
&\hspace*{-,25cm}=&\hspace*{-,25cm}\omega\circ
\pi\bigl(\mathfrak{E}(fagf^*\lambda_{-\gamma})\bigr)\\
&\hspace*{-,25cm}=&\hspace*{-,25cm}\int d\hat{\gamma}\,
\omega\Bigl(\pi\bigl(\hat{\rho}(fagf^*\lambda_{-\gamma})\bigr)\Bigr)=0.
\end{eqnarray*}
 Hence
$\omega\Bigl(\pi\bigl(\mathfrak{E}(fagf^*\lambda_{-\gamma})\bigr)\Bigr)=
0$ for any $\omega\in \B(H)_*$ and
$\pi\bigl(\mathfrak{E}(fagf^*\lambda_{-\gamma})\bigr)=0$. But
$\mathfrak{E}(fagf^*\lambda_{-\gamma})\in A$ and $\pi$ is faithful
on $A$ hence
\begin{equation}\label{eqmor3}\mathfrak{E}(fagf^*\lambda_{-\gamma})=0
 \mbox{ for all }f,g\in L^2(\hat{\Gamma})\cap\C_\infty(\hat{\Gamma}).\end{equation}
We will show that the above equation may be satisfied only if
$a=0$. Let $f_\epsilon\in L^1(\Gamma)$ be an approximation of the
Dirac delta function as used in Theorem 7.8.7 of \cite{Pe}. This
theorem says that for any $y$ of the form $y=fagf^*$ we have the
following norm convergence:
\[\lim_{\varepsilon\rightarrow 0}\int
\mathfrak{E}(y\lambda_{-\gamma})\lambda_{\gamma}\lambda_{f_\epsilon}d\gamma=y.\]
Using \rf{eqmor3} we get
$\mathfrak{E}(y\lambda_{-\gamma})=\mathfrak{E}(fagf^*\lambda_{-\gamma})=0$
hence $fagf^*=0$ for all $f,g\in
L^2(\hat{\Gamma})\cap\C_\infty(\hat{\Gamma})$.
 This immediately implies that $a=0$ and shows that $\pi$ is faithful on
$A^\Psi$. A similar argument shows that faithfulness of $\pi$ on
$A^\Psi$ implies its faithfulness on $A$.
\end{proof}
\begin{defin}\label{quantmor} Let $(A,\rho,\Psi)$, $(A',\rho',\Psi')$
be deformation data with
  groups $\Gamma$ and $\Gamma'$ respectively. Let $\phi:\Gamma\mapsto \Gamma'$ be
  a surjective continuous homomorphism, $\phi^T:\hat{\Gamma}'\mapsto
  \hat{\Gamma}$ the dual homomorphism
 and $\pi\in \Mor(A,A')$. We say that $(\pi,\phi)$ is a
morphism of  deformation data $(A,\rho,\Psi)$ and
$(A',\rho',\Psi')$ if:
\begin{itemize}
\item $\Psi\circ(\phi^T\times \phi^T)=\Psi'$;
\item ${\rho'}_{\phi(\gamma)}\pi(a)=\pi(
\rho_{\gamma}(a))$.
\end{itemize}
\end{defin}
Using universal properties of crossed products, we see that a
morphism $(\pi,\phi)$  of the deformation data induces the
morphism $\pi^\phi\in \Mor(A\rtimes \Gamma;A'\rtimes \Gamma')$ of
crossed products. One can check that
  $\pi^\phi$ satisfies the assumptions of Proposition
\ref{gmor3} with the $\Gamma$-product $(A\rtimes_\rho
\Gamma,\lambda,\hat{\rho})$ and the $\Gamma'$-product
$(A'\rtimes_{\rho'} \Gamma',\lambda',\hat{\rho}')$. This property
is not spoiled by the deformation procedure. Applying Proposition
\ref{gmor3} and Theorem \ref{faithfulrep} to the morphism
$\pi^\phi\in \Mor (A\rtimes_\rho \Gamma;A'\rtimes_{\rho'}
\Gamma')$, $\Gamma$-product $(A\rtimes_\rho
\Gamma,\lambda,\hat{\rho}^\Psi)$ and $\Gamma'$-product
$(A'\rtimes_{\rho'} \Gamma',\lambda',\hat{\rho}'^{\Psi'})$ we get
\begin{stwr}\label{gmor4}  Let  $(\pi,\phi)$ be a
 morphism of deformation data $(A,\rho,\Psi)$ and $(A',\rho',\Psi')$
 and let $\pi^\phi\in \Mor (A\rtimes_\rho\Gamma;A'\rtimes_{\rho'}\Gamma')$
be the induced morphism of the crossed products considered above.
 Then $\pi^\phi(A^\Psi)\subset \M(A'^{\Psi'})$ and
  $\pi^\phi|_{A^\Psi}\in \Mor(A^\Psi;A'^{\Psi'})$.
 Morphism $\pi\in \Mor(A;A')$ is injective if and only if  so is
 $\pi^\phi|_{A^\Psi}\in\Mor(A^\Psi;A'^{\Psi'})$ and
  $\pi(A)=A'$ if and only if $\pi^\phi({A^\Psi})={A'}^{\Psi'}$.
\end{stwr}

 Let $(\mathcal{I},\Gamma,\rho_{\mathcal{I}})$, $(A,\Gamma,\rho)$,
$(A',\Gamma,\rho')$ be dynamical systems and let
\begin{equation}\label{r1} 0\rightarrow \mathcal{I}\rightarrow A
\build{\rightarrow}_{}^{\pi} A' \rightarrow 0\end{equation}
 be an exact sequence of $\C^*$-algebras which is $\Gamma$-equivariant.
 Morphism $\pi$ induces a surjective morphism
 $\pi\in \Mor(A\rtimes_{\rho} \Gamma;A'\rtimes_{\rho'}
 \Gamma)$. It sends $A$ to $A'$ by means of $\pi$ and it is identity on
 $\C^*(\Gamma)$. Its kernel can be identified with
 $\mathcal{I}\rtimes_{\rho_{\mathcal{I}}}
 \Gamma$ so we have an exact sequence of crossed product $\C^*$-algebras:
 \begin{equation} \label{ses} 0\rightarrow \mathcal{I}\rtimes_{\rho_{\mathcal{I}}} \Gamma
 \rightarrow A\rtimes_\rho \Gamma\build{\rightarrow}_{}^{\pi }A'\rtimes_{\rho'} \Gamma\rightarrow
 0.
 \end{equation}
Note that $\pi\in \Mor(A\rtimes_{\rho} \Gamma;A'\rtimes_{\rho'}
\Gamma)$ satisfies the assumptions of Proposition \ref{exct} with
the $\Gamma$-products $(A\rtimes_{\rho}
\Gamma,\lambda,\hat{\rho})$ and $(A'\rtimes_{\rho'}
\Gamma,\lambda',\hat{\rho}')$. This property is not spoiled by the
deformation procedure. Hence applying Proposition $\ref{exct}$ to
the $\Gamma$-products $(A\rtimes_{\rho}
\Gamma,\lambda,\hat{\rho}^\Psi)$ and $(A'\rtimes_{\rho'}
\Gamma,\lambda',\hat{\rho}'^\Psi)$ we obtain
\begin{twr} \label{doklskr} Let $(\mathcal{I},\Gamma,\rho_{\mathcal{I}})$, $(A,\Gamma,\rho)$,
$(A',\Gamma,\rho')$ be dynamical systems. Let  \[0\rightarrow
\mathcal{I}\rightarrow A \build{\rightarrow}_{}^{\pi} A'
\rightarrow 0\] be an exact sequence of $\C^*$-algebras which is
$\Gamma$-equivariant, $\Psi$  a $2$-cocycle on the dual group $
\hat{\Gamma}$ and $\mathcal{I}^{\,\Psi}$, $A^\Psi$, $A'^{\Psi}$
the Landstad algebras constructed  from the deformation data
$(\mathcal{I},\rho_{\mathcal{I}},\Psi)$, $(A,\rho,\Psi)$,
$(A',\rho',\Psi)$. Then we have the $\Gamma$-equivariant exact
sequence:
\[0\rightarrow \mathcal{I}^{\,\Psi}\rightarrow
A^\Psi\build{\longrightarrow}_{}^{\pi^\Psi}A'^{\,\Psi}\rightarrow
0\]
 where the morphism $\pi^\Psi\in \Mor(A^\Psi;A'^\Psi)$ is the
restriction of the morphism $\pi\in \Mor(A\rtimes_\rho
\Gamma,A'\rtimes_{\rho'} \Gamma)$ to the Landstad algebra
$A^\Psi\subset \M(A\rtimes_\rho \Gamma)$.
\end{twr}
\end{subsection}
\begin{subsection}{Preservation of nuclearity.}
\begin{twr} Let $(A,\rho,\Psi)$ be the deformation data which gives rise
to the Landstad algebra $A^\Psi$. $\C^*$-algebra $A$ is nuclear if
and only if $A^\Psi$ is.
\end{twr}
The proof follows from the equality $A\rtimes_\rho
\Gamma=A^\Psi\rtimes_{\rho^\Psi} \Gamma$ (Proposition \ref{eqgam})
and the following:
\begin{twr} Let $A$ be a $\C^*$-algebra with an action $\rho$ of
an abelian group $\Gamma$.  Then $A$ is nuclear if and only if
$A\rtimes_\rho \Gamma$ is nuclear.
\end{twr}
The above theorem can be deduced from Theorem $3.3$
 and Theorem $3.16$ of \cite{take}.
\end{subsection}
\begin{subsection} {$\mathcal{K}$-theory in the case of $\Gamma=\mathbb{R}^n$.}
 In this section
we will prove the invariance of $ \mathcal{K}$-groups under the
Rieffel deformation in the case of $\Gamma=\mathbb{R}^n$. The tool
we use is the analogue of the Thom isomorphism due to Connes
\cite{Co2}:
\begin{twr}\label{thom}\,\, Let $A$ be a $\C^*$-algebra,
and  $\rho$  an action of $\mathbb{R}^n$ on $A$. Then
\[\mathcal{K}_i(A)\simeq\mathcal{K}_{i+n}(A\rtimes_{\rho}\mathbb{R}^n).\]
\end{twr}
\begin{twr}\label{kth} Let $(A,\mathbb{R}^n,\rho)$ be a dynamical
system and let $(A,\rho,\Psi)$ be the deformation data giving rise
to the Landstad algebra $A^\Psi$. Then
\[ \mathcal{K}_i(A)\simeq\mathcal{K}_i(A^\Psi).\]
\end{twr}
\begin{proof}
Proposition \ref{eqgam} asserts that
\[A\rtimes_\rho \mathbb{R}^n\simeq A^\Psi\rtimes_{\rho^\Psi}\mathbb{R}^n.\] Hence
using Theorem \ref{thom} we get
\[\mathcal{K}_i(A)\simeq\mathcal{K}_{i+n}(A\rtimes_\rho\mathbb{R}^n)
\simeq\mathcal{K}_{i+n}(A^\Psi\rtimes_{\rho^\Psi}\mathbb{R}^n)\simeq\mathcal{K}_i(A^\Psi).\]
\end{proof}
\end{subsection}
\end{section}
\begin{section}{Rieffel deformation of locally compact groups.}\label{qg}
\begin{subsection} {From an abelian subgroup with a dual
$2$-cocycle to a quantum group.}\label{fromg} In this section we
shall apply our deformation procedure to the algebra of functions
on a locally compact group $G$. First we shall fix a notation and
introduce auxiliary objects. Let $G\ni g\mapsto
R_g\in\B\bigl(L^2(G)\bigr)$ be the right regular representation of
$G$ on Hilbert space $L^2(G)$ of the  right invariant Haar
measure. Let $\C_\infty(G)\subset \B\bigl(L^2(G)\bigr)$ be the
$\C^*$-algebra of continuous functions on $G$ vanishing at
infinity, $\C^*_r(G)\subset \B\bigl(L^2(G)\bigr)$  the reduced
group $\C^*$-algebra generated by $R_g$ and  $V\in
\B\bigl(L^2(G\times G)\bigr)$ the Kac-Takesaki operator:
$Vf(g,g')=f(gg',g')$  for any $f\in L^2(G\times G)$. By
$\Delta_G\in \Mor\bigl(\C_\infty(G);\C_\infty(G)\otimes
\C_\infty(G)\bigr)$ we will denote the comultiplication on
$\C_\infty(G)$. It is known that the Kac-Takesaki operator $V$ is
an element of $\M\bigl(\C^*_r(G)\otimes\C_\infty(G)\bigr)$ which
implements comultiplication:\[\Delta_G(f)=V(f\otimes 1)V^* \] for
any $f\in\C_\infty(G)$. Let $\Gamma\subset G$ be an abelian
subgroup of $G$, $\hat{\Gamma}$  its dual group and
$\Delta_{\hat\Gamma}\in
\Mor\bigl(\C_\infty(\hat\Gamma);\C_\infty(\hat\Gamma)\otimes
\C_\infty(\hat\Gamma)\bigr)$ the comultiplication on
$\C_\infty(\hat\Gamma)$. Let $ \pi^R\in\Mor\bigr(\C^*(
 \Gamma);\C^*_r(G)\bigr)$ be a morphism induced by the following
 representation of the group $\Gamma$:
\[\Gamma\ni\gamma\mapsto R_\gamma \in \M\bigl(\C^*_r(G)\bigr).\]
 Identifying  $\C^*(\Gamma)$ with
$\C_\infty(\Hat\Gamma)$ we get  $\pi^R\in\Mor\bigr(\C_\infty(
 \Hat\Gamma);\C^*_r(G)\bigr)$.

Let us fix a $2$-cocycle $\Psi$ on the group $\hat{\Gamma}$. Our
objective is to show that an  action of $\Gamma^2$ on the
$\C^*$-algebra $\C_\infty (G)$ given by the left and  right shifts
and a $2$-cocycle on $\hat{\Gamma}^2$ determined by $\Psi$, give
rise to a quantum group. We shall describe this construction step
by step.

Let $\rho^R$ be the action of $\Gamma$ on $\C_\infty(G)$ given by
right shifts: $\rho^R_\gamma(f)(g)=f(g\gamma)$ for any $f\in
 \C_\infty(G)$. Let  $B^R\,$  be the crossed
product $\C^*$- algebra $\C_\infty(G)\rtimes_{\rho^R} \Gamma$ and
$(B^R,\lambda,\hat{\rho})$ the standard $\Gamma$-product structure
on it. The standard embeddings of $\C_\infty(G)$ and
$\C_\infty(\hat\Gamma)$ into $\M(B^R)$ enable us to treat
$(\pi^R\otimes {\rm id})\Psi$ and $V^*(1\otimes f)V$ (where
$f\in\C_\infty(\hat\Gamma)$) as elements of
$\M\bigl(\C^*_r(G)\otimes B^R\bigr)$. One can show that
$V^*(1\otimes \lambda_\gamma)V=R_\gamma\otimes
\lambda_\gamma\mbox{ for all } \gamma\in \Gamma$, which implies
that
\begin{equation} \label{com1}V^*(1\otimes f)V=(\pi^R\otimes {\rm
id})\Delta_{\hat\Gamma}(f)\end{equation} for any  $f\in
\C_\infty(\hat\Gamma)$. Using $\Psi$ we deform the standard
$\Gamma$-product structure on $B^R$ to
$(B^R,\lambda,\hat\rho^\Psi)$.
 \begin{stwr}\label{inv1}Let $(B^R,\lambda,\hat{\rho}^\Psi)$ be the deformed
  $\Gamma$-product and $V (\pi^R\otimes {\rm id})\Psi\in \M(\C^*_r(G)\otimes
B^R)$ the unitary element considered above. Then
 $V (\pi^R\otimes {\rm id})\Psi$
is invariant with respect to the action ${\rm id}
\otimes\hat{\rho}^\Psi$.
\end{stwr} \begin{proof} The $2$-cocycle equation for $\Psi$ implies that:
\begin{eqnarray*}({\rm id}\otimes \hat{\rho}^\Psi_{ \hat{\gamma}} )\Psi
&\hspace*{-,25cm}=&\hspace*{-,25cm}({\rm id}\otimes \hat{\rho}_{
\hat{\gamma}})\Psi\\
&\hspace*{-,25cm}=&\hspace*{-,25cm}(I\otimes U_{
\hat{\gamma}})^*\Delta_{\hat{\Gamma}} ( U_{ \hat{\gamma}})\Psi.
\end{eqnarray*}
The second leg of $V$ is invariant with respect to  the action $
\hat{\rho}$ hence
\begin{eqnarray*}({\rm id}\otimes \hat{\rho}^\Psi_{ \hat{\gamma}} )V
&\hspace*{-,25cm}=&\hspace*{-,25cm}(I\otimes U_{
\hat{\gamma}}^*)(({\rm id}\otimes \hat{\rho}_{ \hat{\gamma}} )V)(
I\otimes U_{ \hat{\gamma}})\\
&\hspace*{-,25cm}=&\hspace*{-,25cm}(I\otimes U_{
\hat{\gamma}}^*)V( I\otimes U_{ \hat{\gamma}}) =V(\pi^R\otimes{\rm
id})\Delta_{\hat{\Gamma}}(U_{ \hat{\gamma}})^*( I\otimes U_{
\hat{\gamma}}).\end{eqnarray*} The last equality follows from
\eqref{com1}. Finally
\begin{eqnarray*} ({\rm id}\otimes \hat{\rho}^\Psi_{ \hat{\gamma}}
)[V(\pi^R\otimes {\rm id})\Psi]&\\
&\hspace*{-4cm}=&\hspace*{-2cm}V(\pi^R\otimes{\rm
id})\Delta_{\hat{\Gamma}} (U_{ \hat{\gamma}})^*( I\otimes U_{
\hat{\gamma}})(I\otimes U_{ \hat{\gamma}})^*(\pi^R\otimes{\rm
id})\Delta_{\hat{\Gamma}}
( U_{ \hat{\gamma}})(\pi^R\otimes {\rm id})\Psi\\
&\hspace*{-4cm}=&\hspace*{-2cm}V(\pi^R\otimes {\rm
id})\Psi\end{eqnarray*} where in the last equality we used the
fact that $ U_{ \hat{\gamma}}$ is unitary.
\end{proof}
Let $\rho^L$ be the action of $\Gamma$ on $\C_\infty(G)$ given by
left shifts: $\rho^L_\gamma(f)(g)=f(\gamma^{-1}g)$ for any $f\in
 \C_\infty(G)$. Let $B^L$ be the crossed
product $\C^*$- algebra $\C_\infty(G)\rtimes_{\rho^L} \Gamma$ and
let $(B^L,\lambda,\hat{\rho})$ be the standard $\Gamma$-product
structure on it. For any
$\hat{\gamma}_1,\hat{\gamma}_2\in\hat{\Gamma}$ we set
\[\Psi^\star(\hat{\gamma}_1,\hat{\gamma}_2)=
\overline{\Psi(\hat{\gamma}_1,-\hat{\gamma}_1-\hat{\gamma}_2)}.\]
This defines a function $\Psi^\star\in \C_b(\hat{\Gamma}^2)$.  The
standard embeddings of $\C_\infty(G)$ and
$\C_\infty(\hat{\Gamma})$ into $\M(B^L)$ enable us to treat
$(\pi^R\otimes {\rm id})\Psi^\star V$  and  $V(1\otimes f)V^*$
(where $f\in\C_\infty(\hat\Gamma)$)
  as elements of
$\M\bigl(\C^*_r(G)\otimes B^R\bigr)$. One can show that
$V(1\otimes \lambda_\gamma)V^*=R_\gamma\otimes \lambda_\gamma$ for
all $\gamma\in \Gamma$, which implies that
\begin{equation}\label{com11} V(1\otimes f)V^*=(\pi^R\otimes {\rm id})
\Delta_{\hat{\Gamma}}(f)\end{equation} for any $f\in
\C_\infty(\hat{\Gamma})$. Let $\widetilde{\Psi}$ denote a
$2$-cocycle  defined by the formula:
\[ \widetilde{\Psi}(\hat{\gamma}_1, \hat{\gamma}_2)\equiv
\overline{\Psi(-\hat{\gamma}_1, -\hat{\gamma}_2)}\] for any
$\hat{\gamma}_1,\hat{\gamma}_2\in \hat{\Gamma}$. Using
$\widetilde{\Psi}$ we deform the standard $\Gamma$-product
structure on $B^L$ to $\left(B^L,\lambda,
\hat{\rho}^{\widetilde{\Psi}}\right)$.
\begin{stwr}\label{inv2} Let $\left(B^L,\lambda,
\hat{\rho}^{\widetilde{\Psi}}\right)$ be the deformed
$\Gamma$-product and
 $(\pi^R\otimes {\rm id})\Psi^\star V\in \M(\C^*_r(G)\otimes B^L)$
 the unitary element considered above. Then $(\pi^R\otimes {\rm
id})\Psi^\star V\in \M(\C^*_r(G)\otimes B^L)$ is invariant with
respect to the action ${\rm id}\otimes
\hat{\rho}^{\widetilde{\Psi}}$.
\end{stwr}\begin{proof} One can check that
\[\Psi^\star(\hat{\gamma}_1,\hat{\gamma}_2+\hat{\gamma})=
\overline{\Psi(\hat{\gamma}_1,-\hat{\gamma}_1-\hat{\gamma}_2-\hat{\gamma})}
=\widetilde{\Psi}(\hat{\gamma}_2,\hat{\gamma})
\overline{\widetilde{\Psi}(\hat{\gamma}_1+\hat{\gamma}_2,\hat{\gamma})}\Psi^\star(\hat{\gamma}_1,\hat{\gamma}_2).\]
Hence
\[({\rm id}\otimes \hat{\rho}^{\widetilde{\Psi}}_{\hat{\gamma}})\Psi^\star
 =(I\otimes\widetilde{U}_{
\hat{\gamma}})\Delta_{\hat{\Gamma}}(\widetilde{U}_{
\hat{\gamma}})^*\Psi^\star.
\] Moreover \[({\rm id}\otimes
\hat{\rho}^{\widetilde{\Psi}}_{\hat{\gamma}})V=(I\otimes
\widetilde{U}_{ \hat{\gamma}})^*V(I\otimes \widetilde{U}_{
\hat{\gamma}})=(I\otimes \widetilde{U}_{
\hat{\gamma}})^*(\pi^R\otimes{\rm id})
\Delta_{\hat{\Gamma}}(\widetilde{U}_{ \hat{\gamma}})V.\] Following
the proof of Proposition \ref{inv1} we get our assertion.
\end{proof}

Let $\rho$ denote the action of $\Gamma^2$ on $ \C_\infty(G)$
given by the left and right shifts, $B$ the crossed product
$\C^*$-algebra $ \C_\infty(G)\rtimes_\rho \Gamma^2$ and
$(B,\lambda,\hat{\rho})$ the standard  $\Gamma^2$-product. The
standard embedding of $\C_\infty(G)$ into $\M(B)$ applied  to the
second leg of $V\in \M\bigl(\C_r^*(G)\otimes \C_\infty(G)\bigr)$
 embeds  $V$  into $\M\bigl(\C^*_r(G)\otimes
B\bigr)$. We have two embeddings $\lambda^L$ and $\lambda^R$ of $
\C_\infty( \hat{\Gamma})$ into $\M(B)$ corresponding to the left
and the right action of $\Gamma$. Moreover by equations
\eqref{com1} and \eqref{com11} we have:
\begin{equation}\label{impl0}\begin{array}{c} V(1\otimes
\lambda^L(f))V^*=(\pi^R\otimes {\rm
\lambda^L})\Delta_{\hat\Gamma}(f) \\V^*(1\otimes
\lambda^R(f))V=(\pi^R\otimes {\rm
\lambda^R})\Delta_{\hat\Gamma}(f)\end{array}\end{equation} for any
$f\in \C_\infty(\hat{\Gamma})$. Note also that:
\begin{equation}\label{impl1}({ \rm id}\otimes\lambda_{\gamma_{1},\gamma_{2}})
V({\rm
id}\otimes\lambda^*_{\gamma_{1},\gamma_{2}})=(R_{-\gamma_{1}}\otimes
I)V(R_{\gamma_{2}}\otimes I).\end{equation} Let us introduce
elements $\Psi^L$ and $ \Psi^R$:
\begin{equation}\label{psilr}\Psi^L=(\pi^R\otimes \lambda^L)(\Psi^\star),\,\,
\Psi^R=(\pi^R\otimes \lambda^R)(\Psi)\,\in\,
\M\bigl(\C^*_r(G)\otimes B\bigr).\end{equation} Multiplying
$\Psi^L,V $ and $ \Psi^R $ we get the unitary  element:
\begin{equation}\label{vaupsi} V^\Psi=\Psi^LV \Psi^R\in \M\bigl(\C^*_r(G)\otimes
B\bigr).\end{equation} Using the $2$-cocycle $
\widetilde{\Psi}\otimes \Psi$ on $\hat{\Gamma}^2$ we deform the
standard $\Gamma^2$-product structure on $B$ to
$(B,\lambda,\hat{\rho}^{ \widetilde{\Psi}\otimes \Psi})$.
\begin{stwr} \label{wninv}Let $(B,\lambda, \hat{\rho}^{\widetilde{\Psi}\otimes \Psi})$
be the deformed $\Gamma^2$-product structure and  $V^\Psi\in
\M\bigl(\C^*_r(G)\otimes B\bigr)$ the unitary element given by
\eqref{vaupsi}. Then $V^\Psi$ is invariant with respect to the
action $ {\rm id}\otimes \hat{\rho}^{ \widetilde{\Psi}\otimes
\Psi}$. Moreover, for any $\gamma_1,\gamma_2\in \Gamma$ we have
\begin{equation} \label{impl2}\rm(id\otimes\lambda_{\gamma_{1},\gamma_{2}})
V^\Psi\rm(id\otimes\lambda^*_{\gamma_{1},\gamma_{2}})=(R_{-\gamma_{1}}\otimes
I)V^\Psi(R_{\gamma_{2}}\otimes I).\end{equation}
\end{stwr}
\begin{proof} Invariance of $V^\Psi$ with respect to
the action $ {\rm id}\otimes \hat{\rho}^{ \widetilde{\Psi}\otimes
\Psi}$ follows easily from Propositions \ref{inv1} and \ref{inv2}.
The group $\Gamma$ is abelian, hence \begin{eqnarray*} ({\rm
id}\otimes\lambda_{\gamma_{1},\gamma_{2}})V^\Psi({\rm
id}\otimes\lambda^*_{\gamma_{1},\gamma_{2}})
&\hspace*{-,25cm}=&\hspace*{-,25cm}({\rm
id}\otimes\lambda_{\gamma_{1},\gamma_{2}})\Psi^LV \Psi^R({\rm
id}\otimes\lambda^*_{\gamma_{1},\gamma_{2}})\\
&\hspace*{-.25cm}=&\hspace*{-.25cm}\Psi^L({\rm
id}\otimes\lambda_{\gamma_{1},\gamma_{2}})V ({\rm
id}\otimes\lambda^*_{\gamma_{1},\gamma_{2}})\Psi^R.
\end{eqnarray*} Using \eqref{impl1} we get \begin{eqnarray*}
({\rm id}\otimes\lambda_{\gamma_{1},\gamma_{2}})V^\Psi({\rm id}\otimes\lambda^*_{\gamma_{1},\gamma_{2}})
&\hspace*{-.25cm}=&\hspace*{-.25cm}\Psi^L(R_{-\gamma_{1}}\otimes
I)V(R_{\gamma_{2}}\otimes I)\Psi^R\\
&\hspace*{-.25cm}=&\hspace*{-.25cm}(R_{-\gamma_{1}}\otimes
I)\Psi^LV\Psi^R(R_{\gamma_{2}}\otimes I)\\
&\hspace*{-.25cm}=&\hspace*{-.25cm}(R_{-\gamma_{1}}\otimes
I)V^\Psi(R_{\gamma_{2}}\otimes I).
\end{eqnarray*} This proves \eqref{impl2}.
\end{proof}
 The first leg of $V^\Psi$ belongs to $\C^*_r(G)$ so
it acts on $L^2(G)$. It is well-known that slices of  Kac-Takesaki
operator $V$ by normal functionals $\omega\in\B(L^2(G))_*$ give a
dense subspace of $ \C_\infty(G)$ (see \cite{BS}, Section 2). We
will show that slices of $V^\Psi$ give a dense subspace  of $
\C_\infty(G)^{\widetilde{\Psi}\otimes \Psi}$.
\begin{twr}\label{lc} Let $(B,\lambda, \hat{\rho}^{\widetilde{\Psi}\otimes \Psi})$
be the deformed $\Gamma^2$-product structure and  $V^\Psi\in
\M\bigl(\C^*_r(G)\otimes B\bigr)$ the unitary operator given by
\eqref{vaupsi}. Then
\[\mathcal{V}=\bigl\{(\omega\otimes \id)V^\Psi:\omega\in
\B\bigl(L^2(G)\bigr)_*\bigr\}\] is a norm dense subset of $
\C_\infty(G)^{\widetilde{\Psi}\otimes \Psi}$.
\end{twr}
\begin{proof} We need to check that for any $\omega\in
\B\bigl(L^2(G)\bigr)_*$
 the element  $(\omega\otimes \id)V^\Psi\in\M(B)$
satisfies Landstad conditions for $\Gamma^2$-product
$(B,\lambda,\hat{\rho}^{ \widetilde{\Psi}\otimes \Psi})$. The
first Landstad condition is equivalent to the invariance of the
second leg of $V^\Psi$ with respect to the action $\hat{\rho}^{
\widetilde{\Psi}\otimes \Psi}$ (Proposition \ref{wninv}). Using
\eqref{impl2} we get
\begin{equation}\label{roinv}\lambda_{\gamma_{1},\gamma_{2}}[(\omega\otimes \rm
id)V^\Psi]\lambda_{\gamma_{1},\gamma_{2}}^*=(R_{\gamma_2}\cdot\omega\cdot
R_{-\gamma_1}\otimes \id)V^\Psi
\end{equation} for any
$\gamma_1,\gamma_2\in \Gamma$. The norm continuity of the map
\[\Gamma^2\ni(\gamma_1,\gamma_2)\mapsto
R_{\gamma_2}\cdot\omega\cdot R_{-\gamma_1}\in \B(L^2(G))_*\]
implies that  $(\omega\otimes \rm id)V^\Psi$ satisfies  the second
Landstad condition. To check the third Landstad condition we need
to show that
 \begin{equation}\label{ll3}f_1[(\omega\otimes {\rm
id})V^\Psi]f_2\in
 B
 \end{equation} for any  $f_1,f_2\in \C_{\infty}(\hat{\Gamma}\times
\hat{\Gamma})$. Let us consider the set
\begin{equation}\label{denscr}\mathcal{W}=\{f_1[(\omega\otimes {\rm id})V^\Psi]f_2:f_1,f_2\in
\C_{\infty}(\hat{\Gamma}\times \hat{\Gamma}),
\omega\in\B\bigl(L^2(G)\bigr)_*\}^{\rm \rm cls}.\end{equation} We
will prove that $\mathcal{W}=B$ which is a stronger property than
\eqref{ll3}. Taking for $\omega\in\B\bigl(L^2(G)\bigr)_*$ elements
of the form $ \pi^R(h_3)\cdot\mu\cdot\pi^R(h_4)$, for
$f_1\in\C_{\infty}(\hat{\Gamma}\times \hat{\Gamma})$
 elements $\lambda^R(h_1)\lambda^L(h_2)$ where $h_1,h_2\in\C_{\infty}(\hat{\Gamma})$
  and  similarly for $f_2$ we do not change
 the closed linear span. Thus we have:
\begin{eqnarray*}\mathcal{W}=\{\lambda^R(h_1)\lambda^L(h_2)[((\pi^R(h_3)\cdot\mu\cdot\pi^R(h_4))\otimes
{\rm id})(V^\Psi)]\lambda^R(h_5)\lambda^L(h_6):&\\&\hspace*{-5,5cm}
h_1,h_2\ldots,h_6\in \C_\infty(\hat{\Gamma}),\mu\in
\B\bigl(L^2(G)\bigr)_*\}^{\rm cls}.
\end{eqnarray*}
Note that
\[\lambda^R(h_1)\lambda^L(h_2)[((\pi^R(h_3)\cdot\mu\cdot\pi^R(h_4))\otimes
{\rm id})(V^\Psi)
]\lambda^R(h_5)\lambda^L(h_6)\hspace*{3cm}\]\[=\lambda^R(h_1)[(\mu\otimes
{\rm id})(\pi^R\otimes \lambda^L)(\Psi ^\star(h_4\otimes h_2))
V(\pi^R\otimes \lambda^R)(\Psi (h_3\otimes h_5))]\lambda^L(h_6)\]
hence $\mathcal{W}$ coincides with the following set:
\begin{eqnarray*}\{\lambda^R(h_1)[(\mu\otimes {\rm
id})(\pi^R\otimes \lambda^L)(\Psi ^\star(h_4\otimes h_2))
V(\pi^R\otimes \lambda^R)(\Psi (h_3\otimes
h_5))]\lambda^L(h_6):\\&\hspace{-7cm}h_1,h_2,\ldots,h_6\in
\C_\infty(\hat{\Gamma}),\mu\in \B\bigl(L^2(G)\bigr)_*\}^{\rm cls}.
\end{eqnarray*}
Using the fact that $\Psi$ and $\Psi^{\star}$ are unitary we get
\begin{eqnarray*}\mathcal{W}&\hspace*{-0,25cm}=&\hspace*{-0,25cm}\{\lambda^R(h_1)[(\mu\otimes {\rm
id})(\pi^R\otimes \lambda^L)(h_4\otimes h_2) V(\pi^R\otimes
\lambda^R) (h_3\otimes
h_5)]\lambda^L(h_6):\\&\hspace*{-0,25cm}&\hspace*{4cm}h_1,h_2\ldots,h_6\in
\C_\infty(\hat{\Gamma}),\mu\in \B\bigl(L^2(G)\bigr)_*\}^{\rm
cls}\\&\hspace*{-0,25cm}=&\hspace*{-0,25cm}\{\lambda^R(h_1)\lambda^L(h_2)[((\pi^R(h_3)\cdot\mu\cdot\pi^R(h_4))\otimes
{\rm id})(V)
]\lambda^R(h_5)\lambda^L(h_6):\\&\hspace*{-0,25cm}&\hspace*{4cm}h_1,h_2\ldots,h_6\in
\C_\infty(\hat{\Gamma}),\mu\in \B\bigl(L^2(G)\bigr)_*\}^{\rm cls}
\end{eqnarray*}
Now again
\begin{eqnarray*}\{\lambda^R(h_1)\lambda^L(h_2)[((\pi^R(h_3)\cdot\mu\cdot\pi^R(h_4))\otimes
{\rm id})(V)
]\lambda^R(h_5)\lambda^L(h_6):&\\&\hspace*{-5,5cm}h_1,h_2\ldots,h_6\in
\C_\infty(\hat{\Gamma}),\mu\in \B\bigl(L^2(G)\bigr)_*\}^{\rm
cls}&\\&\hspace*{-8.3cm}=\{f_1[(\omega\otimes {\rm
id})V]f_2:f_1,f_2\in \C_{\infty}(\hat{\Gamma}\times
\hat{\Gamma}),\omega\in\B\bigl(L^2(G)\bigr)_*\}^{\rm
cls}\end{eqnarray*} hence we get
\begin{eqnarray*}\mathcal{W}=\{f_1[(\omega\otimes {\rm
id})V]f_2:f_1,f_2\in \C_{\infty}(\hat{\Gamma}\times
\hat{\Gamma}),\omega\in\B\bigl(L^2(G)\bigr)_*\}^{\rm cls}.
\end{eqnarray*}
The set $\{(\omega\otimes {\rm id})V:\omega\in
\B\bigl(L^2(G)\bigr)_*\}$ is dense in $\C_\infty(G)$ which shows
that $\mathcal{W}=B$ and  proves formula \rf{denscr}.

We see that the elements of the set $\mathcal{V}$ satisfy the
Landstad conditions. To prove that $ \mathcal{V}$ is dense in $
\C_\infty(G)^{\widetilde{\Psi}\otimes \Psi}$
 we use Lemma \ref{dsw}. According to \eqref{roinv},
 $\mathcal{V}$ is a $\rho^{\widetilde{\Psi}\otimes
\Psi}$-invariant subspace of $
\C_\infty(G)^{\widetilde{\Psi}\otimes \Psi}$. Moreover we have
that $\left(\C^*(\Gamma^2) \mathcal{V}\C^*(\Gamma^2)\right)^{\rm
cls}=\mathcal{W}=B$. Hence the assumptions of Lemma \ref{dsw} are
satisfied and we get the required density.
\end{proof}
\begin{uw}\label{kan}
The representation of $\C_\infty(G)$ on $L^2(G)$ is covariant. The
action of $\Gamma^2$ is implemented by the left and
 right shifts:  $L_{\gamma_1},R_{\gamma_2}\in
\B\bigl(L^2(G)\bigr)$, where by $L_g\in\B(L^2(G))$ we  understand
the unitarized left shift. More precisely, let
$\delta:G\rightarrow\mathbb{R}^+$ be the modular function for the
right Haar measure. Then $L_g\in\B(L^2(G))$ is a unitary given by:
\[(L_gf)(g')=\delta(g)^{\frac{1}{2}}f(g^{-1}g')\] for any $g,g'\in G$ and $f\in L^2(G)$.
This covariant representation of $\C_\infty(G)$ induces the
representation of crossed product $B=
\C_\infty(G)\rtimes_\rho\Gamma^2$, which we denote by $\pi^{\rm
can}$. Clearly it is faithful on $\C_\infty(G)$,  hence by Theorem
\ref{faithfulrep} it is faithful on
$\C_\infty(G)^{\widetilde{\Psi}\otimes \Psi}$. \end{uw} Let us
introduce the unitary operator: \begin{equation}\label{utw}
W=({\rm id}\otimes \pi^{\rm can})V^\Psi\in \B\bigl(L^2(G)\otimes
L^2(G)\bigr).\end{equation}
\begin{twr}\label{multun} The unitary operator $W\in\B\bigl(L^2(G)\otimes L^2(G)\bigr)$
 considered above  satisfies the pentagonal
equation: \[W_{12}^*W_{23}W_{12}=W_{13}W_{23}.\]
\end{twr}
\begin{uw} A similar  construction of the operator $W$ and the proof that it satisfies
the pentagonal equation was given by Enock-Vainerman in \cite{V}
and independently by Landstad in \cite{lan1}. We included the
following proof for the completeness of the exposition.
\end{uw}
\begin{proof}
Let us introduce two unitary operators $X,Y\in\B(L^2(G)\otimes
L^2(G))$:
\begin{equation}\label{XY0} X=({\rm id}\otimes \pi^{\rm
can})(\Psi^R)\,\,\,, \,\,\,Y=({\rm id}\otimes \pi^{\rm
can})(\Psi^L)\end{equation} where $\Psi^R,\Psi^L\in
\M(\C^*_r(G)\otimes B)$ are elements defined by \eqref{psilr}. Note
that
\begin{equation}\label{xynal}
X\in\M(\C_r^*(G)\otimes\C_r^*(G)),\,\,
Y\in\M(\C_r^*(G)\otimes\C_l^*(G)),\end{equation} hence
$W=YVX\in\M(\C_r^*(G)\otimes\mathcal{K})$ where $\mathcal{K}$ is
the algebra of compact operators acting on $L^2(G)$. Inserting
$\hat{\gamma}_3\rightarrow(-\hat{\gamma}_1-\hat{\gamma}_2-\hat{\gamma}_3)$
into the $2$-cocycle condition
\begin{equation}\label{eqpsi}\Psi(\hat{\gamma}_1,\hat{\gamma}_2+\hat{\gamma}_3)\Psi(\hat{\gamma}_2,\hat{\gamma}_3)
=\Psi(\hat{\gamma}_1+\hat{\gamma}_2,\hat{\gamma}_3)\Psi(
\hat{\gamma}_1,\hat{\gamma}_2) \end{equation} and taking the
complex conjugate we get
\[
\overline{\Psi(\hat{\gamma}_1,-\hat{\gamma}_1-\hat{\gamma}_3)}
\overline{\Psi(\hat{\gamma}_2,-\hat{\gamma}_1-\hat{\gamma}_3-\hat{\gamma}_2)}=
\overline{\Psi(\hat{\gamma}_1,\hat{\gamma}_2)}
\overline{\Psi(\hat{\gamma}_1+\hat{\gamma}_2,-\hat{\gamma}_1-\hat{\gamma}_2-\hat{\gamma}_3)}.\]
This  implies that
\begin{equation}\label{eqpsitil}\overline{\Psi}(\hat{\gamma}_1,\hat{\gamma}_2)\Psi^\star(\hat{\gamma}_1+\hat{\gamma}_2,\hat{\gamma}_3)=
\Psi^\star(\hat{\gamma}_1,\hat{\gamma}_3)\Psi^\star(\hat{\gamma}_2,\hat{\gamma}_1+\hat{\gamma}_3)
\end{equation}
where
$\Psi^\star(\hat{\gamma}_1,\hat{\gamma}_2)=\overline{\Psi(\hat{\gamma}_1,-\hat{\gamma}_1-\hat{\gamma}_2)}$.
Using equations \eqref{XY0}, \eqref{eqpsi}, \eqref{eqpsitil} and the
fact that $V$ implements the coproduct we obtain:
\begin{eqnarray*}
X_{12}^{*}V_{12}^{*}Y_{23}V_{12}&\hspace*{-0,2cm}=&\hspace*{-0,1cm}Y_{13}V_{13}Y_{23}V_{13}^*\\
V^*_{12}X_{23}V_{12}X_{12}&\hspace*{-0,2cm}=&\hspace*{-0,1cm}V^*_{23}X_{13}V_{23}X_{23}
.\end{eqnarray*} Now we can check the pentagonal equation:
\begin{eqnarray*} W_{12}^{*}W_{23}W_{12}
&\hspace*{-.25cm}=&\hspace*{-.25cm}
X_{12}^{*}V_{12}^{*}Y_{12}^{*}W_{23}Y_{12}V_{12}X_{12}\\
&\hspace*{-.25cm}=&\hspace*{-.25cm}X_{12}^{*}V_{12}^{*}Y_{23}V_{23}X_{23}V_{12}X_{12}\\
&\hspace*{-.25cm}=&\hspace*{-.25cm}(X_{12}^{*}V_{12}^{*}Y_{23}V_{12})
(V_{12}^*V_{23}V_{12})(V_{12}^*X_{23}V_{12}X_{12})\\
&\hspace*{-.25cm}=&\hspace*{-.25cm}
(Y_{13}V_{13}Y_{23}V_{13}^*)(V_{13}V_{23})(V^*_{23}X_{13}V_{23}X_{23})\\
&\hspace*{-.25cm}=&\hspace*{-.25cm}Y_{13}V_{13}Y_{23}X_{13}V_{23}X_{23}=
(Y_{13}V_{13}X_{13})(Y_{23}V_{23}X_{23})
=W_{13}W_{23}.\end{eqnarray*} In the second equality we used the
fact that the second leg of element $Y$ commutes with the first leg
of $W$ (see \eqref{xynal}).
\end{proof}

Our next aim is to show that $W$ is manageable.  For all
$\hat{\gamma}\in \hat{\Gamma}$ we set
$u(\hat{\gamma})=\Psi(-\hat{\gamma},\hat{\gamma})$. It defines a
function $u\in\C_b(\hat{\Gamma})$. Applying $\pi^R\in
\Mor\bigl(\C_\infty ( \hat{\Gamma});\C^*_r(G)\bigr)$ to $u\in
\M\bigl(\C_\infty ( \hat{\Gamma})\bigr)$ we get the unitary
operator:
\begin{equation}\label{uop}J=\pi^R(u)\in \M(\C^*_r(G))\subset
\B\bigl(L^2(G)\bigr).\end{equation}
\begin{twr} \label{twr4}
 Let $W\in \B\bigl(L^2(G)\otimes L^2(G)\bigr)$ be the  multiplicative
 unitary and  $J\in \B\bigl(L^2(G)\bigr)$ be the unitary
 operator \eqref{uop}. Then $W$ is manageable.
Operators $Q$ and $\widetilde{W}$ entering the Definition 1.2 of
\cite{W5} equal respectively  $1$ and
 $(J\otimes 1)W^*(J^*\otimes 1)$.
 \end{twr}
\begin{uw}  The presented proof seems to be simpler than
the  Landstad's proof given in \cite{lan1}. In what follows we
shall use the bracket notation for the scalar product: let $H$ be
a Hilbert space,  $x,y\in H$, and $T\in\B(H)$. Then  $(x|T|y)$
denotes the scalar product $(x|Ty)$.
\end{uw}
\begin{proof}
Let $x,y,z,t\in L^{2}(G)$,
$\gamma_{1},\gamma_{2},\gamma_{3},\gamma_{4}\in \Gamma$. The
Kac-Takesaki operator is manageable, therefore
\begin{eqnarray*}(x\otimes t|(R_{\gamma_{1}}\otimes
L_{\gamma_{2}})V|(R_{\gamma_{3}}\otimes R_{\gamma_{4}})|z\otimes
y)\\&\hspace*{-4cm}=&\hspace*{-2cm}(R_{-\gamma_{1}}x\otimes
L_{-\gamma_{2}}t|V|R_{\gamma_{3}}z\otimes R_{\gamma_{4}}
y)\\&\hspace*{-4cm}=&\hspace*{-2cm}(\overline{R_{\gamma_{3}}z}\otimes
L_{-\gamma_{2}}t|V^{*}|\overline{R_{-\gamma_{1}}x}\otimes
R_{\gamma_{4}}
y)\\&\hspace*{-4cm}=&\hspace*{-2cm}(\overline{z}\otimes
t|(R_{-\gamma_{3}}\otimes
L_{\gamma_{2}})V^{*}(R_{-\gamma_{1}}\otimes
R_{\gamma_{4}})|\overline{x}\otimes y).\end{eqnarray*} Using well
known equalities
\[V^{*}(I\otimes R_{g})V=R_{g}\otimes R_{g}\]
\[V(I\otimes L_{g})V^{*}=R_{g}\otimes L_{g}\]
and commutativity of $\Gamma$ we get the following formula:
\begin{equation}
\begin{array}{ccc} (\overline{z}\otimes t|(R_{-\gamma_{3}}\otimes
L_{\gamma_{2}})V^{*}(R_{-\gamma_{1}}\otimes
R_{\gamma_{4}})|\overline{x}\otimes
y)&&\\&\hspace*{-10cm}=&\hspace*{-5cm}(\overline{z}\otimes
t|(R_{-\gamma_{3}+\gamma_{4}}\otimes
R_{\gamma_{4}})V^{*}(R_{-\gamma_{1}+\gamma_{2}}\otimes
L_{\gamma_{2}})|\overline{x}\otimes y).
\end{array}
\end{equation}
Hence:
\begin{equation}\label{48}
\begin{array}{ccc} (x\otimes t|(R_{\gamma_{1}}\otimes
L_{\gamma_{2}})V(R_{\gamma_{3}}\otimes R_{\gamma_{4}})|z\otimes
y)&&\\&\hspace*{-10cm}=&\hspace*{-5cm} (\overline{z}\otimes
t|(R_{-\gamma_{3}+\gamma_{4}}\otimes
R_{\gamma_{4}})V^{*}(R_{-\gamma_{1}+\gamma_{2}}\otimes
L_{\gamma_{2}})|\overline{x}\otimes y).
\end{array}\end{equation}
Using continuity arguments, this equality will be extended. We
will repeatedly use the identifications
$\C^*(\Gamma^2)=\C_\infty(\Hat\Gamma^2)=\C_\infty(\Hat\Gamma)\otimes\C_\infty(\Hat\Gamma)$
etc. Let $u_\gamma $ be a unitary generator of $\C^*(\Gamma)$. Let
us define the following morphisms:
\begin{eqnarray*}\Phi_{1}^{R}\in \Mor(\C^*(\Gamma)\otimes
\C^*(\Gamma);\C^{*}_{r}(G)\otimes \C^{*}_{r}(G)):&&
\Phi_{1}^{R}(u_{\gamma_{1}}\otimes
u_{\gamma_{2}})=R_{\gamma_{1}}\otimes R_{\gamma_{2}},\\
\Phi_{1}^{L}\in \Mor(\C^*(\Gamma)\otimes
\C^*(\Gamma);\C^{*}_{r}(G)\otimes \C^{*}_{l}(G)):&&
\Phi_{1}^{L}(u_{\gamma_{1}}\otimes
u_{\gamma_{2}})=R_{\gamma_{1}}\otimes
L_{\gamma_{2}}\end{eqnarray*} and  automorphism $ \Theta\in
{\rm Aut}(
 \C_{\infty}( \hat{\Gamma}^2))$ given by the formula:
\[\Theta(f)(
\hat{\gamma}_1,\hat{\gamma}_2)=f(-\hat{\gamma}_1,\hat{\gamma}_1+\hat{\gamma}_2)\]
for any $f\in  \C_{\infty}( \hat{\Gamma}^2)$. One can check that
\[\Theta(u_{\gamma_1}\otimes
u_{\gamma_2})=u_{-\gamma_1+\gamma_2}\otimes u_{\gamma_2}.\] Using
the above morphisms we reformulate
\eqref{48}:\begin{eqnarray*}(x\otimes
t|(\Phi_{1}^{L}(u_{\gamma_{1}}\otimes
u_{\gamma_{2}})V\Phi_{1}^{R}(u_{\gamma_{3}}\otimes
u_{\gamma_{4}})|z\otimes
y)\\&\hspace*{-8cm}=&\hspace*{-4cm}(\overline{z}\otimes
t|\Phi_{1}^{R}\circ\Theta(u_{\gamma_{3}}\otimes
u_{\gamma_{4}})V^{*}\Phi_{l}^{L}\circ\Theta(u_{\gamma_{1}}\otimes
u_{\gamma_{2}})|\overline{x}\otimes y).\end{eqnarray*} By
linearity and continuity we get
\begin{eqnarray*}(x\otimes
t|\Phi_{1}^{L}(f)V\Phi_{1}^{R}(g)|z\otimes
y)\\&\hspace*{-6cm}=&\hspace*{-3cm}(\overline{z}\otimes
t|\Phi_{1}^{R}\circ\Theta(g)V^{*}\Phi_{1}^{L}\circ\Theta(f)|\overline{x}\otimes
y)\end{eqnarray*} for any $f,g\in \M(\C_\infty(\Hat\Gamma)\otimes
\C_\infty(\Hat\Gamma))$. In particular
\begin{eqnarray*}(x\otimes
t|\Phi_{1}^{L}(\Psi^\star)V\Phi_{1}^{R}(\Psi)|z\otimes
y)\\&\hspace*{-6cm}=&\hspace*{-3cm}(\overline{z}\otimes
t|\Phi_{1}^{R}\circ\Theta(\Psi)V^{*}\Phi_{1}^{L}\circ\Theta(\Psi^\star)|\overline{x}\otimes
y)\end{eqnarray*} It is easy to see that  $X=\Phi_{1}^{R}(\Psi)$,
$Y=\Phi_{1}^{L}(\Psi^\star)$ and
 $\Theta(\Psi)=\overline{\Psi}(u\otimes
I)$ where $X$ and $Y$ are given by \eqref{XY0}. Therefore
\begin{equation} \label{tild1}\Phi_{1}^{R}\circ\Theta(\Psi)=\Phi_{1}^{R}\left(\overline{\Psi}(u\otimes I
)\right)=( J\otimes I)X^*.\end{equation}
 Similarly we prove that
\begin{equation}\label{tild2}\Phi_{1}^{L}\circ\Theta(\Psi^\star)=Y^*(J^*\otimes
I)\end{equation} and  finally we get
\begin{eqnarray*}(x\otimes
t|\Phi_{1}^{L}(\Psi^\star)V\Phi_{1}^{R}(\Psi)|z\otimes
y)\\&\hspace*{-3cm}=&\hspace*{-1,5cm}(\overline{z}\otimes t|(
J\otimes I)X^*V^*Y^*(J^*\otimes I)|\overline{x}\otimes
y).\end{eqnarray*} This shows that \[\widetilde{W}=\left((J\otimes
I)YVX( J^*\otimes I)\right)^*=(J\otimes I)W^*(J^*\otimes
I)\mbox{\,\,\, and \,\,\,}Q=1.
\]
\end{proof}
\begin{stwr}\label{man} Let $W\in\B(L^2(G)\otimes L^2(G))$ and $J\in\B(L^2(G))$ be the
 unitaries  defined in \eqref{utw} and \eqref{uop} respectively.
 Let $x,y$ be vectors in $L^2(G)$ and $\omega_{x,y}\in
\B(L^2(G))_*$ a functional given by $\omega_{x,y}(T)=\left(
x|T|y\right)$ for any $T\in\B(L^2(G))$. Then we have
\begin{equation} \label{menst}[(\omega_{x,y}\otimes {\rm id})W]^*=
(\omega_{J^*\bar{x},J^*\bar{y}}\otimes {\rm
id})(W).\end{equation}\end{stwr}
\begin{proof} Using  manageability of $W$ we get:
\begin{eqnarray*}(\omega_{x,y}\otimes {\rm
id})W&\hspace*{-0,25cm}=&\hspace*{-0,25cm}(\omega_{\bar{y},\bar{x}}\otimes
{\rm id})(\widetilde{W})=(\omega_{\bar{y},\bar{x}}\otimes {\rm
id})((J\otimes 1)W^*(J^*\otimes
1))\\&\hspace*{-0,25cm}=&\hspace*{-0,25cm}(\omega_{J^*\bar{y},J^*\bar{x}}\otimes
{\rm id})(W^*)=[(\omega_{J^*\bar{x},J^*\bar{y}}\otimes {\rm
id})(W)]^*.\end{eqnarray*}\end{proof} Let $A$ be a $\C^*$-algebra
obtained by slicing the first leg of a manageable multiplicative
unitary $W\in \B\bigl(L^2(G)\otimes L^2(G)\bigr)$:
\[A=\overline{\bigl\{(\omega\otimes\rm{id})W:\omega\in\B(L^2(G))_*\bigr\}}^{\|\cdot\|}.\]
Theorem 1.5 of \cite{W5} shows that $A$ carry the structure of a
quantum group. The comultiplication on $A$ is given by the
formula:
\[A\ni a\mapsto W(a\otimes I)W^*\in \M(A\otimes A).\]
At the same time, using the morphism $\pi^{\rm
can}\in\Rep(B;L^2(G))$ introduced in Remark \ref{kan} we can
faithfully represent $\C_\infty(G)^{\widetilde{\Psi}\otimes \Psi}$
on $L^2(G)$. By Theorem \ref{lc}  $\pi^{\rm
can}(\C_\infty(G)^{\widetilde{\Psi}\otimes \Psi})=A$, hence we can
transport the structure of a quantum group from $A$ to
$\C_\infty(G)^{\widetilde{\Psi}\otimes \Psi}$. Our next objective
is to present a  useful formula for comultiplication on $
\C_\infty(G)^{\widetilde{\Psi}\otimes \Psi}$ which does not use
multiplicative unitary $W$. The construction is done in two steps.
\begin{itemize}
  \item Let $\rho$ be the action of\,\, $\Gamma^2$ on $ \C_{\infty}(G)$ given by left and right shifts along the subgroup
  $\Gamma\subset G$. The comultiplication is covariant:
  \[\Delta_G\bigl(\rho_{\gamma_1,\gamma_2}(f)\bigr)=
  (\rho_{\gamma_1,0}\otimes \rho_{0,\gamma_2})\bigl(\Delta_G(f)\bigr)\] for any $f\in
  \C_\infty(G)$. Therefore, it induces a morphism of crossed products:
  \[\Delta\in {\rm Mor}\bigl( \C_{\infty}(G)\rtimes \Gamma^2;
  \C_{\infty}(G)\rtimes \Gamma^2\otimes  \C_{\infty}(G)\rtimes
  \Gamma^2\bigr).\]
  $\Delta$ restricted to $\C_{\infty}(G)\subset\M\bigl( \C_{\infty}(G)\rtimes
    \Gamma^2)$ coincides with  $\Delta_G$ and  $\Delta$
   restricted to
 $\C_\infty( \hat{\Gamma}^2)\subset\M\bigl( \C_{\infty}(G)\rtimes
    \Gamma^2)$ is given by \[\Delta(h)=(\lambda^L\otimes\lambda^R)h\in
   \M\bigl(  \C_{\infty}(G)\rtimes \Gamma^2\otimes  \C_{\infty}(G)\rtimes
\Gamma^2\bigr)\] where
$\lambda^L,\lambda^R\in\Mor(\C_\infty(\Hat\Gamma),\C_\infty(G)\rtimes\Gamma^2)$
are morphisms introduced after the proof of Proposition \ref{inv2}
and $h\in\C_\infty( \hat{\Gamma}^2)$.
  \item Let $\Psi$ be a $2$-cocycle on $ \hat{\Gamma}$.
  Recall that  $\Psi^{\star}\in  \M\bigl( \C_\infty( \hat{\Gamma}^2)\bigr)$
is defined by
  \[\Psi^{\star}(\hat{\gamma}_1,\hat{\gamma}_2)=\overline{\Psi(\hat{\gamma}_1,-\hat{\gamma}_1-\hat{\gamma}_2)}.\]
Let us introduce the unitary element $\Upsilon\in\M\bigl(
\C_{\infty}(G)\rtimes
  \Gamma^2\otimes \C_{\infty}(G)\rtimes
  \Gamma^2\bigr)$:
\[ \Upsilon=(\lambda^R\otimes \lambda^L)\Psi^\star\]
and a morphism $ \Delta^\Psi\in\Mor\bigl( \C_{\infty}(G)\rtimes
\Gamma^2;   \C_{\infty}(G)\rtimes \Gamma^2\otimes
\C_{\infty}(G)\rtimes \Gamma^2\bigr)$ given by the formula
\begin{equation}\label{delpsidef}\Delta^\Psi(a)=\Upsilon
\Delta(a)\Upsilon^* \end{equation} for any
$a\in\C_{\infty}(G)\rtimes \Gamma^2$.
\end{itemize}

\begin{twr}\label{forcom} Let $\Delta^\Psi\in\Mor\bigl(
\C_{\infty}(G)\rtimes \Gamma^2;   \C_{\infty}(G)\rtimes
\Gamma^2\otimes  \C_{\infty}(G)\rtimes \Gamma^2\bigr)$ be the
morphism defined by formula \eqref{delpsidef}. For all $a\in
\C_\infty(G)^{\widetilde{\Psi}\otimes \Psi}$ we have
\[\Delta^\Psi(a)\in \M\bigl( \C_\infty(G)^{\widetilde{\Psi}\otimes
\Psi}\otimes
 \C_\infty(G)^{\widetilde{\Psi}\otimes \Psi}\bigr)\] and
\[\Delta^\Psi|_{ \C_\infty(G)^{\widetilde{\Psi}\otimes \Psi}}\in
\Mor\bigl( \C_\infty(G)^{\widetilde{\Psi}\otimes \Psi};
\C_\infty(G)^{\widetilde{\Psi}\otimes \Psi}\otimes
 \C_\infty(G)^{\widetilde{\Psi}\otimes \Psi}\bigr).\] Moreover
 $\Delta^\Psi|_{ \C_\infty(G)^{\widetilde{\Psi}\otimes \Psi}}$
coincides with the comultiplication implemented by $W$:
\[\C_\infty(G)^{\widetilde{\Psi}\otimes \Psi}\ni a\mapsto W(a\otimes 1)W^*
\in \M\bigl(\C_\infty(G)^{\widetilde{\Psi}\otimes
\Psi}\otimes\C_\infty(G)^{\widetilde{\Psi}\otimes \Psi}\bigr) .\]
\end{twr}
\begin{proof} By Theorem 1.5 of \cite{W5} it is enough
to show that  \[({\rm id} \otimes
\Delta^\Psi)V^\Psi=V^\Psi_{12}V^\Psi_{13}.\] From the definition
of $ \Delta$ it follows that \begin{eqnarray*} ({\rm id} \otimes
\Delta)V^\Psi&\hspace*{-.25cm}=&\hspace*{-.25cm}({\rm id} \otimes
\Delta)\bigl((\pi^R\otimes \lambda^L)(\Psi^{\star})V(\pi^R\otimes
\lambda^R)(\Psi)\bigr)\\&\hspace*{-.25cm}=&\hspace*{-.25cm}\bigl((\pi^R\otimes
\lambda^L)\Psi^{\star}V\bigr)_{12}\bigl(V(\pi^R\otimes
\lambda^R)\Psi\bigr)_{13} .\end{eqnarray*} Hence
\begin{eqnarray*}({\rm id} \otimes
\Delta^\Psi)V^\Psi =(1\otimes\Upsilon)\bigl((\pi^R\otimes
\lambda^L)\Psi^{\star}V\bigr)_{12}\bigl(V((\pi^R\otimes
\lambda^R)\Psi\bigr)_{13})(1\otimes\Upsilon^*) .\end{eqnarray*} By
equation \eqref{impl0} we get
\begin{eqnarray*}\bigl(1\otimes\Upsilon\bigr)V_{12}&\hspace*{-.25cm}=&\hspace*{-.25cm}
V_{12}\Bigl((\pi^R\otimes\lambda^R\otimes\lambda^L)\circ(\Delta_{\hat{\Gamma}}\otimes
{\rm id})\bigl( \Psi^\star\bigr)\Bigr)
\end{eqnarray*} and
\begin{eqnarray*}V_{13}\bigl(1\otimes\Upsilon\bigr)&\hspace*{-.25cm}=&\hspace*{-.25cm}
\Bigl((\pi^R\otimes\lambda^R\otimes\lambda^L)\circ(\sigma\otimes{\rm
id})\circ ({\rm id}\otimes\Delta_{\hat{\Gamma}} )
\bigl(\Psi^\star\bigr)\Bigr)V_{13}
\end{eqnarray*} where $\sigma$ is the flip operator.
Therefore
\begin{eqnarray*}({\rm id} \otimes
\Delta^\Psi)V^\Psi
&\hspace*{-.25cm}=&\hspace*{-.25cm}\bigl((\pi^R\otimes
\lambda^L)\Psi^{\star}V\bigr)_{12}\Bigl((\pi^R\otimes\lambda^R\otimes\lambda^L)\circ(\Delta_{\hat{\Gamma}}\otimes
{\rm id})\bigl(
\Psi^\star\bigr)\Bigr)\\&\hspace*{-.25cm}\times&\hspace*{-.25cm}
\Bigl((\pi^R\otimes\lambda^R\otimes\lambda^L)\circ(\sigma\otimes{\rm
id})\circ ({\rm id}\otimes\Delta_{\hat{\Gamma}} )
\bigl(\Psi^\star\bigr)^*\Bigr)\bigl(V(\pi^R\otimes
\lambda^R)\Psi\bigr)_{13} .\end{eqnarray*} We compute
\begin{eqnarray*}
\overline{\Psi(\hat{\gamma}_1+\hat{\gamma}_2,-\hat{\gamma}_1-\hat{\gamma}_2-\hat{\gamma}_3)}
\Psi(\hat{\gamma}_2,-\hat{\gamma}_1-\hat{\gamma}_2-\hat{\gamma}_3)
\\&\hspace*{-12cm}=&\hspace*{-6cm}\Psi(\hat{\gamma}_1,\hat{\gamma}_2)
\overline{\Psi(\hat{\gamma}_1,-\hat{\gamma}_1-\hat{\gamma}_3)}
\overline{\Psi(\hat{\gamma}_2,-\hat{\gamma}_1-\hat{\gamma}_2-\hat{\gamma}_3)}
\Psi(\hat{\gamma}_2,-\hat{\gamma}_1-\hat{\gamma}_2-\hat{\gamma}_3)
\\&\hspace*{-10cm}&\hspace*{0cm}=\Psi(\hat{\gamma}_1,\hat{\gamma}_2)
\overline{\Psi(\hat{\gamma}_1,-\hat{\gamma}_1-\hat{\gamma}_3)}.\end{eqnarray*}
The above equality implies that
\begin{eqnarray*}
(\pi^R\otimes\lambda^R\otimes\lambda^L)\circ(\Delta_{\hat{\Gamma}}\otimes
{\rm id})\bigl(
\Psi^\star\bigr)\\&\hspace*{-6cm}\times&\hspace*{-3cm}
(\pi^R\otimes\lambda^R\otimes\lambda^L)\circ(\sigma\otimes{\rm
id})\circ ({\rm id}\otimes\Delta_{\hat{\Gamma}} )
\bigl(\Psi^\star\bigr)^*\\&\hspace*{-2cm}=&\hspace*{-1cm}\bigl((\pi^R\otimes
\lambda^R)\Psi\bigr)_{12}\bigl((\pi^R\otimes
\lambda^L)\Psi^\star\bigr)_{13}.\end{eqnarray*} Hence
 \begin{eqnarray*}({\rm id} \otimes\Delta^\Psi)V^\Psi
&\hspace*{-.25cm}=&\hspace*{-.25cm}\bigl((\pi^R\otimes
\lambda^L)\Psi^{\star}V(\pi^R\otimes
\lambda^R)\Psi\bigr)_{12}\bigl((\pi^R\otimes \lambda^L)\Psi^\star
V(\pi^R\otimes
\lambda^R)\Psi\bigr)_{13}\\&\hspace*{-.25cm}=&\hspace*{-.25cm}V^\Psi_{12}V^\Psi_{13}.
\end{eqnarray*} This ends the proof.
\end{proof}
\end{subsection}
\begin{subsection}{Dual quantum group.}
 Let $G$ be a locally compact group, $\Gamma$ an abelian subgroup of $G$ and
 $\Psi$ a $2$-cocycle on $\hat{\Gamma}$. Using the results of
 previous sections we can construct the quantum group
 $\bigl( \C_\infty(G)^{\widetilde{\Psi}\otimes \Psi},\Delta^\Psi\bigr)$
 and the multiplicative unitary $W\in \B\bigl(L^2(G)\otimes L^2(G)\bigr)$.
  In this section we will investigate the dual quantum group in the sense of duality
given by $W$. Our objective is to show that this is the twist, in
the sense of M. Enock and L. Vainerman (see \cite{V}), of the
canonical quantum group structure on the reduced group
$\C^*$-algebra $\C^*_r(G)$.
\begin{twr}\label{dual}
Let $W\in \B\bigl(L^2(G)\otimes L^2(G)\bigl)$ be a manageable
multiplicative unitary \eqref{utw} and $(\Hat{A},\Hat\Delta_{\Hat
A})$ a quantum group obtained by slicing the second leg of $W$:
\[\Hat{A}=\overline{\bigl\{(\id\otimes\omega)(W^*):\omega\in
\B\bigl(L^2(G)\bigr)\bigr\}}^{\|\cdot\|}.\] Then \begin{itemize}
\item[1.]
$\hat{A}=\C_r^*(G)$. \item[2.] The comultiplication on $\hat{A}$
is given by
\[\hat{A}\ni a\mapsto \hat{\Delta}_{\hat{A}}(a)=
\Sigma X^*\Sigma\hat{\Delta}(a)\Sigma X \Sigma\in
\M(\hat{A}\otimes \hat{A})\] where $\hat{\Delta}$ is the canonical
comultiplication on $\C^*_r(G)$ and $X$ is given by \eqref{XY0}.
\item[3.] The coinverse on $\hat{A}$ is given by
\[\hat{\kappa}_{\hat{A}}(a)=J\hat{\kappa}(a)J^*\]
where $\hat{\kappa}$ is the canonical coinverse on $\C^*_r(G)$ and
$J$ is given by \eqref{uop}.
\end{itemize}
\end{twr}
The proof was communicated to the author by  S.L. Woronowicz.
\begin{proof}
Using equation \eqref{impl2} we get
\begin{equation} \label{strcont}R_{\gamma_1}[(\id\otimes \omega)W]R_{\gamma_2}=({\rm id}\otimes
R_{-\gamma_2}L_{\gamma_1}\cdot\omega\cdot
L_{-\gamma_1}R_{\gamma_2})W\end{equation} for any
$\gamma_1,\gamma_2\in\Gamma$. Therefore $R_\gamma\in\B(L^2(G))$ is
a multiplier of $\hat{A}$ and representation:
 \[\Gamma\ni\gamma\mapsto R_\gamma\in \M(\hat{A})\] is
strictly continuous. This representation induces a morphism which we
denote by
 $\chi\in \rm Mor \bigl(\C_\infty(\hat{\Gamma}), \hat{A}\bigr)$.
Applying it to $\Psi$ and  $\Psi^\star$ we obtain
\begin{eqnarray*}
X&\hspace*{-0,25cm}=&\hspace*{-0,25cm}(\chi\otimes \chi)(\Psi) \in
\M( \hat{A}\otimes
\hat{A})\\Y&\hspace*{-0,25cm}=&\hspace*{-0,25cm}(\chi\otimes
\pi^L)(\Psi^\star) \in  \M( \hat{A}\otimes \mathcal{K}).
\end{eqnarray*} Recall that $W\in\M(\hat{A}\otimes A)$, hence
\begin{equation}\label{pent} V=Y^*WX^*\in \M( \hat{A}\otimes
\mathcal{K})\end{equation} which immediately implies that
$V_{12}^*V_{23},V_{12}V_{23}^*\in \M\bigl(\hat{A}\otimes \mathcal{K}
\otimes
 \C_{ \infty}(G)\bigr)$. The pentagonal equation  for $V$ together with
 \eqref{pent} gives
\[ V_{13}=V_{12}^*V_{23}V_{12}V_{23}^*\in \M\bigl(\hat{A}\otimes
\mathcal{K} \otimes
 \C_{ \infty}(G)\bigr),\] therefore \begin{equation}\label{vgen} V\in
\M\bigl(\hat{A} \otimes  \C_{ \infty}(G)\bigr).\end{equation}
Similarly we prove that
\begin{equation}\label{wgen}W\in \M\bigl(\C_{r}^*(G)\otimes A\bigr).\end{equation}

Formula \eqref{vgen} and point 6 of Theorem 1.6 of \cite{W5} imply
that the natural representation of $\C_r^*(G)$ on $L^2(G)$ is in
fact an element of $ \Mor\bigl(\C_r^*(G),\hat{A}\bigr)$.
Similarly, \rf{wgen} implies that the natural representation of
$\hat{A}$ on $L^2(G)$ is an element of
$\Mor\bigl(\hat{A},\C_r^*(G)\bigr)$. The general properties of
morphisms gives \begin{eqnarray*}
\overline{\C_r^*(G)\hat{A}}&\hspace*{-0,25cm}=&\hspace*{-0,25cm}\hat{A}\\
\overline{\hat{A}\C_r^*(G)}&\hspace*{-0,25cm}=&\hspace*{-0,25cm}\C_r^*(G).\end{eqnarray*}
But $\C_r^*(G)$ and $\hat{A}$ are closed under the star operation,
hence $\hat{A}=\hat{A}^*=\overline{\hat{A}\C^*_r(G)}=\C_r^*(G)$,
which proves point 1 of our theorem. To prove point 2 we recall
that the  comultiplication on $\hat{A}$ is
 implemented by $\Sigma W^*\Sigma$, hence
\begin{eqnarray*}\hat{\Delta}_{\hat{A}}(a)&\hspace*{-.25cm}=&\hspace*{-.25cm}
\Sigma X^*V^*Y^*(I\otimes a)YVX \Sigma\\
&\hspace*{-.25cm}=&\hspace*{-.25cm}\Sigma X^*V^*(I\otimes a)VX\Sigma
\\&\hspace*{-.25cm}=&\hspace*{-.25cm}(\Sigma X^*\Sigma)\hat{\Delta}(a)(\Sigma X\Sigma).
\end{eqnarray*}
Point three follows from Proposition \ref{man}.\end{proof}
\end{subsection}
\begin{subsection}{Haar measure.} Let $G$ be a locally compact group,
$\Gamma$ an  abelian subgroup of $G$ and $\Psi$ a $2$-cocycle on
$\hat {\Gamma}$. Throughout this section we shall assume that the
modular function $\delta$ on $G$ restricted to $\Gamma$ is
identically equal to $1$. Let $\bigl(
\C_\infty(G)^{\widetilde{\Psi}\otimes \Psi},\Delta^\Psi\bigr)$ be
the quantum group that we considered previously. In what follows
we will identify $\C_\infty(G)^{\widetilde{\Psi}\otimes \Psi}$
with its image in $B\bigl(L^2(G)\bigr)$.
\begin{defin}\label{qmap} Let $f\in  \C_\infty(G)$ and $R_g\in \B\bigl(L^2(G)\bigr)$
be the right regular representation of  group $G$. We say that $f$
is quantizable if there exists $\omega\in \B\bigl(L^2(G)\bigr)_*$
such that $f(g)=\omega(R_g)$ for any $g\in G$. Given a quantizable
function $f$ we introduce an operator
 $\mathcal{Q}(f)\in \C_\infty(G)^{\widetilde{\Psi}\otimes \Psi}\subset \B(L^2(G))$ given by:
\[\mathcal{Q}(f)=(\omega\otimes {\rm id})W\in \C_\infty(G)^{\widetilde{\Psi}\otimes \Psi}.\]
\end{defin}
Note that the equation $f(g)=\omega(R_g)$ does not determine
$\omega\in\B(L^2(G))_*$. Nevertheless, the operator
$\mathcal{Q}(f)$ does not depend on the choice of the functional
that gives rise to $f$. It is easy to see that the vector space of
quantizable  functions equipped with the pointwise multiplication
forms an algebra which in the literature is called the {\it
Fourier Algebra}. We use the term {\it quantizable function} to
stress that with such an $f$ we can associate the operator
$\mathcal{Q}(f)=(\omega\otimes\id)W\in\B(L^2(G))$.
\begin{twr}\label{quant} Let $( \C_\infty(G)^{\widetilde{\Psi}\otimes
\Psi},\Delta^\Psi)$ be the quantum group with multiplicative
unitary $W\in \B\bigl(L^2(G)\otimes L^2(G)\bigr)$ considered
above.
 Let $f,h\in \C_\infty(G)$ be quantizable functions given by functionals
$\omega\in \B\bigl(L^2(G)\bigr)_*$ and $\mu\in \B(L^2(G))_*$
respectively. They yield operators
$\mathcal{Q}(f),\mathcal{Q}(h)\in \B\bigl(L^2(G)\bigr)$. Assume
that $h\in L^2(G)$. Then $\mathcal{Q}(f)h\in L^2(G)$ is a
quantizable function and
\begin{equation}\label{beauty}
\mathcal{Q}\bigl(\mathcal{Q}(f)h\bigr)=\mathcal{Q}(f)\mathcal{Q}(h).\end{equation}
\end{twr} \begin{proof}
Note that
\begin{eqnarray*} \mathcal{Q}(f)\mathcal{Q}(h)&\hspace*{-.25cm}=&\hspace*{-.25cm}(\omega\otimes
\rm{id})(W)(\mu\otimes
\rm{id})W\\
&\hspace*{-.25cm}=&\hspace*{-.25cm}(\omega\otimes\mu\otimes
\rm{id})(W_{13}W_{23})
\\&\hspace*{-.25cm}=&\hspace*{-.25cm}(\omega\otimes\mu\otimes
\rm{id})(W_{12}^{*}W_{23}W_{12}).
\end{eqnarray*}
The above calculation shows that $\mathcal{Q}(f)\mathcal{Q}(h)$ is
given by the quantization of the function $k\in  \C_\infty(G)$:
\[k(g)=(\omega\otimes\mu)\bigl(W^{*}(I\otimes
R_{g})W\bigr).\] Using the identity  $W^{*}(I\otimes
R_{g})W=X^{*}(R_{g}\otimes R_{g})X$ we get
\[k(g)=(\omega\otimes\mu)\bigl(X^{*}(R_{g}\otimes R_{g})X\bigr).\]
Therefore, to prove formula \rf{beauty} we need to show that
\begin{equation}\label{need}(\omega\otimes\mu)\bigl(X^{*}(R_{g}\otimes
R_{g})X\bigr)=[\mathcal{Q}(f)h](g).\end{equation} In order to do
that we compute
\begin{eqnarray*} (\omega\otimes
{\rm id})\bigl((R_{\gamma_{1}}\otimes
L_{\gamma_{2}})V(R_{\gamma_{3}}\otimes
R_{\gamma_{4}})\bigr)h(g)&\hspace*{-.25cm}=&\hspace*{-.25cm}L_{\gamma_{2}}f(\gamma_{1}\cdot
\gamma_{3})R_{\gamma_{4}}h(g)\\&\hspace*{-.25cm}=&\hspace*{-.25cm}\delta^{\frac{1}{2}}(\gamma_2)\,f((\gamma_{1}-\gamma_{2})g\gamma_{3})h((-\gamma_{2})g\gamma_{4}).
\end{eqnarray*}
Using the assumption that $\delta(\gamma)=1$ for any
$\gamma\in\Gamma$ we get
\begin{eqnarray*} (\omega\otimes
{\rm id})\bigl((R_{\gamma_{1}}\otimes
L_{\gamma_{2}})V(R_{\gamma_{3}}\otimes
R_{\gamma_{4}})\bigr)h(g)&\hspace*{-.25cm}=&\hspace*{-.25cm}f((\gamma_{1}-\gamma_{2})g\gamma_{3})h((-\gamma_{2})g\gamma_{4}).
\end{eqnarray*}
 The equality $h(g)=\mu(R_{g})$ implies that
\begin{equation}\label{wzor1}
\begin{array}{l}
 \Bigl((\omega\otimes {\rm id})\bigl((R_{\gamma_{1}}\otimes
L_{\gamma_{2}})V(R_{\gamma_{3}}\otimes
R_{\gamma_{4}})\bigr)h\Bigr)(g)
\\\hspace*{3cm}=(\omega\otimes \mu)\bigl((R_{\gamma_{1}-\gamma_{2}}\otimes
R_{-\gamma_{2}})(R_{g}\otimes R_{g})(R_{\gamma_{3}}\otimes
R_{\gamma_{4}})\bigr).
\end{array}
\end{equation}
 Let
$\vartheta\in \rm Aut(\C_{\infty}(\hat{\Gamma}\times
\hat{\Gamma}))$ be the  automorphism given by
\[\vartheta(f)(\hat{ \gamma}_{1},\hat{\gamma}_{2})=
f(\hat{\gamma}_{1},-\hat{\gamma}_{1}-\hat{\gamma}_{2}) \mbox{ for all } f\in\C_{\infty}(\hat{\Gamma}\times
\hat{\Gamma}).\] By continuity, \eqref{wzor1} extends to
\begin{eqnarray*}[(\omega\otimes {\rm id})\bigl((\pi^{R}\otimes
\pi^{L})(f_1)V(\pi^{R}\otimes
\pi^{R})(f_2)\bigr)]h(g)\\&\hspace*{-8cm}=&\hspace*{-4cm}(\omega\otimes
\mu)\bigl((\pi^{R}\otimes \pi^{R})\vartheta(f_1)(R_{g}\otimes
R_{g})(\pi^{R}\otimes \pi^{R})(f_2)\bigr)\end{eqnarray*} for any
$f_1,f_2\in \C_{b}(\hat{\Gamma}\otimes\hat{\Gamma})$. Taking
$f_1=\Psi^\star$ and $f_2=\Psi$ we obtain
$\vartheta(\Psi^\star)=\overline{\Psi}$ and
\[\rm [(\omega\otimes{\rm id})(YVX)]h(g)=(\omega\otimes
\mu)(X^{*}(R_{g}\otimes R_{g})X)\] where $X$ and $Y$ were
introduced in \eqref{XY0}. Recall that $W=YVX$, hence
\[\mathcal{Q}(f)h(g)=(\omega\otimes
\mu)(X^{*}(R_{g}\otimes R_{g})X)=k(g).\] This proves formula
\eqref{need} and ends the proof of our theorem. \end{proof} Let
$f\in\C_\infty(G)$ be a quantizable function i.e.
$f(g)=\omega(R_g)$ for some $\omega\in\B\bigl(L^2(G)\bigr)_*$.
Suppose that $\mathcal{Q}(f)=0$. This means that $ (\omega\otimes
{\rm id})W=0 $ which together with Theorem \ref{dual} shows  that
$\omega(R_g)=0$ for all $g\in G$. Hence $f(g)=0$ for any $g\in G$,
which shows that the quantization map $ \mathcal{Q}$ is injective
and its inverse is well defined. We shall show that the closure of
this inverse  is the GNS map for a Haar measure of $(
\C_\infty(G)^{\widetilde{\Psi}\otimes \Psi},\Delta^\Psi)$. Let us
introduce
 $\mathfrak{N}_0\subset
\C_\infty(G)^{\widetilde{\Psi}\otimes\Psi}$:
\[\mathfrak{N}_0=\left\{\mathcal{Q}(f): f-\mbox{quantizable and } f\in L^2(G)\right\}.\] For
all $\mathcal{Q}(f)\in \mathfrak{N}_0$
 we set
$\eta_0(\mathcal{Q}(f))=f$. This defines a map
$\eta_0:\mathfrak{N}_0 \rightarrow L^2(G)$.
\begin{stwr} Let  $\eta_0$ be the map defined above.
Then this is a densely defined, closable map from $
\C_\infty(G)^{\widetilde{\Psi}\otimes\Psi}$ to $L^2(G)$.
\end{stwr}
\begin{proof}
 Let $ \Psi^\Sigma$ be a $2$-cocycle obtained from $\Psi$ by a
flip of variables: $\Psi^\Sigma(
\hat{\gamma}_1,\hat{\gamma}_2)=\Psi(
\hat{\gamma}_2,\hat{\gamma}_1)$. Let $ \mathcal{Q}^\Sigma$ be the
quantization map related to $\Psi^\Sigma$ . Using the equality
\[(\omega\otimes \mu)\bigl(X^*(R_g\otimes R_g)X\bigr)=(\mu\otimes
\omega)\bigl(\Sigma X^*\Sigma(R_g\otimes R_g)\Sigma
X\Sigma\bigr)\] and Theorem \ref{quant} we see  that for
quantizable, square integrable functions $h,h'\in L^2(G)$ we have
\begin{equation}\label{lok1} \mathcal{Q}(h)h'=\mathcal{Q}^\Sigma(h')h.\end{equation}
Let us assume  that
$\displaystyle\lim_{n\rightarrow\infty}\mathcal{Q}(f_n)=0$ and
$\displaystyle\lim_{n\rightarrow\infty}\eta_0(\mathcal{Q}(f_n))=f$.
Using equation \eqref{lok1} we see that:
\begin{equation}\label{nicearg}\mathcal{Q}^\Sigma(h)f=\lim_{n\rightarrow\infty}\mathcal{Q}^\Sigma(h)f_n=
\lim_{n\rightarrow\infty}\mathcal{Q}(f_n)h= 0\end{equation}
 for all quantizable functions $h\in L^2(G)$.
 To conclude that $f$ is $0$  we have to show that the set of operators
\[\{\mathcal{Q}^\Sigma(h):h\mbox{ is quantizable
and }h\in L^2(G)\}\subset\B(L^2(G))\] separates elements of
$L^2(G)$. In order to do that we introduce a multiplicative
unitary $W^\Sigma$ related to the $2$-cocycle $\Psi^\Sigma$. The
$\C^*$-algebra obtained by the slices of the first leg of
$W^\Sigma$ will be denoted by $A^{\Psi^\Sigma}$. By point 1 of
Theorem 1.5 of \cite{W5},
 $A^{\Psi^\Sigma}$ separates elements of $L^2(G)$, hence it is enough to note that:
\begin{eqnarray*}A^{\Psi^\Sigma}&\hspace*{-0,25cm}=&\hspace*{-0,25cm}\left\{(\omega\otimes{\rm
id})W^\Sigma:\omega\in\B(L^2(G))\right\}^{\rm cls}\\
&\hspace*{-0,25cm}=&\hspace*{-0,25cm}\left\{(\omega_{x,y}\otimes{\rm
id})W^\Sigma:x,y\mbox{  are of compact support }\right\}^{\rm
cls}\\
&\hspace*{-0,25cm}\subset&\hspace*{-0,25cm}\{\mathcal{Q}^\Sigma(h):h\mbox{
is quantizable and }h\in L^2(G)\}^{\rm cls}.
\end{eqnarray*}
The last inclusion follows from the fact that, when $x$ and $y$
are of compact support, then the function $f$ defined by
$f(g)=\omega_{x,y}(R_g)$ is also of compact support.
\end{proof}
The closure of the map $\eta_0$ will be denoted by $\eta$ and its
domain will be denoted by $\mathfrak{N}$.
\begin{stwr}
Let $ \eta:\mathfrak{N}\mapsto L^2(G)$ be the map introduced above.
Then $\mathfrak{N}$ is a left ideal in $
\C_\infty(G)^{\widetilde{\Psi}\otimes\Psi}$ and
$\eta(ab)=a\eta(b)$ for all $a\in
\C_\infty(G)^{\widetilde{\Psi}\otimes\Psi}$ and $b\in \mathfrak{N}.$
\end{stwr}
\begin{proof} Let $b\in\mathfrak{N}$ and $a\in  \C_\infty(G)^{\widetilde{\Psi}\otimes\Psi}$.
Let us fix a sequence of quantizable functions $f_n$ such that
$\displaystyle a=\lim_{n\rightarrow\infty}\mathcal{Q}(f_n)$. Map
$\eta$ is the closure of $\eta_0$, therefore there exists a
sequence $h_m\in \C_\infty(G)$ of quantizable functions such that:
\[ \lim_{m\rightarrow\infty}\mathcal{Q}(h_m)= b\,\,\,\,{\rm
and}\,\,\,\,
\lim_{m\rightarrow\infty}\eta(\mathcal{Q}(h_m))=\eta(b).\] Using
Theorem \ref{quant} we get
$\mathcal{Q}(f_n)\eta(\mathcal{Q}(h_m))=\eta(\mathcal{Q}(f_n)\mathcal{Q}(h_m))$
and
\[a\eta( \mathcal{Q}(h_m))=\lim_{n\rightarrow\infty}\mathcal{Q}(f_n)\eta(\mathcal{Q}(h_m))
=\lim_{n\rightarrow\infty}\eta(\mathcal{Q}(f_n)\mathcal{Q}(h_m)).\]
The closedness of the map $\eta$ implies that
\[a\mathcal{Q}(h_m)\in \mathfrak{N} \mbox{ and
 } \eta(a\mathcal{Q}(h_m))=a\eta(\mathcal{Q}(h_m)).\]
Taking limits with respect to $m$ and using the closedness of
$\eta$ once again we conclude that $ab\in\mathfrak{N}$ and
$\eta(ab)=a\eta(b)$.
\end{proof}
The above proposition shows that the map $
\eta:\mathfrak{N}\mapsto L^2(G)$ is a GNS map. To show that $\eta$
corresponds  to the Haar measure of
$(\C_\infty(G)^{\tilde\Psi\otimes\Psi}, \Delta^\Psi)$ we shall
need the following
\begin{stwr}\label{haar1}
Let $ \eta:\mathfrak{N}\mapsto L^2(G)$ be the map introduced
above. For $a\in \mathfrak{N}$ and $\varphi\in
(\C_\infty(G)^{\widetilde{\Psi}\otimes\Psi})^{*}$ let us consider
their convolution $\varphi*a=({\rm id}\otimes
\varphi)\Delta(a)\in\C_\infty(G)^{\widetilde{\Psi}\otimes\Psi}$.
Then $\varphi*a$ is an element of $ \mathfrak{N}$ and
\begin{equation} \label{rown}  \eta(\varphi*a)=[({\rm id}\otimes
\varphi)W]\eta(a).\end{equation}
\end{stwr} \begin{proof} Recall that with any normal functional $\omega\in
\B\bigl(L^2(G)\bigr)_*$ we can associate a function $f_\omega\in
\C_\infty(G)$ where $f_\omega(g)=\omega(R_g)$. Assume that
$a=\mathcal{Q}(f_\omega)$ for some $f_\omega\in L^2(G)$. In
particular $a\in \mathfrak{N}$ and $\eta(a)=f_\omega$. We compute
\begin{equation}\label{splot}\begin{array}{rcl}\varphi*a&\hspace*{-.25cm}=&\hspace*{-.25cm}
(\id\otimes \varphi)\bigl(W^*(a\otimes 1)W\bigr)\\
&\hspace*{-.25cm}=&\hspace*{-.25cm}({\rm id}\otimes
\varphi)\bigl(W^*((\omega\otimes {\rm id})(W)\otimes 1)W\bigr)
\\ &\hspace*{-.25cm}=&\hspace*{-.25cm}
 (\omega\otimes {\rm id} \otimes \varphi)(W^*_{23}W_{12}W_{23})\\
&\hspace*{-.25cm}=&\hspace*{-.25cm}(\omega\otimes {\rm id}\otimes
\varphi)(W_{12}W_{13})\\
&\hspace*{-.25cm}=&\hspace*{-.25cm}(b\,\cdot\omega\otimes\id)W
\end{array}\end{equation}
where $b=({\rm id}\otimes \varphi)W\in \M(\C^*_r(G))$. Therefore to
prove that $\varphi*a$ is an element of $\mathfrak{N}$ it is enough
to show that $f_{b\,\cdot\omega}\in L^2(G)$ for all
$b\in\M(\C_r^*(G))$. First we check it for $b=R_{g}$. Note that
\begin{eqnarray*}
f_{b\cdot\omega}(g')
&\hspace*{-.25cm}=&\hspace*{-.25cm}b\cdot\omega(R_g)=\omega(R_g'R_{g})
=\omega(R_{g'g})\\&\hspace*{-.25cm}=&\hspace*{-.25cm}f_{\omega}(g'g)=
(R_gf_\omega)(g')=(bf_\omega)(g')\end{eqnarray*} for any $g,g'\in
G$, therefore $f_{b\cdot\omega}=bf_\omega$. By linearity this
equality is satisfied for any $b\in \mbox{span}\{R_g:g\in G\}$. We
 extend it using a continuity argument. There exists a net of
operators $b_i\in \mbox{lin-span}\{R_g:g\in G\}$ strongly
convergent to $b\in \M(\C^*_r(G))$. Functional $\omega$ is
strongly continuous hence
$\displaystyle\lim_ib_i\cdot\omega=b\,\cdot\omega$ in the norm
sense. Therefore $\displaystyle\lim_i
f_{b_i\cdot\omega}=f_{b\,\cdot\omega}$ where $\displaystyle\lim_i$
is taken in the uniform sense. At the same time $\displaystyle
\lim_ib _if_\omega =b f_\omega \mbox{ in the } L^2-\mbox{norm}$,
hence
\[f_{b\,\cdot\omega}(g)=\lim_if_{b _i\cdot\omega }(g)=\lim_ib
_if_\omega (g)=bf_\omega(g)\]  for almost all $g\in G$. This shows
that $f_{b\,\cdot\omega}\in L^2(G)$ and
\begin{equation}\label{fact}f_{b\,\cdot\omega}=b f_\omega
\end{equation} for any $b\in\M(\C_r^*(G))$.
 Using \eqref{splot} and \eqref{fact} we get the following
sequence of equalities:
\begin{eqnarray*}\eta(\varphi*a)&\hspace*{-.25cm}=&\hspace*{-.25cm}
f_{b\,\cdot\omega}=b
f_\omega\\&\hspace*{-.25cm}=&\hspace*{-.25cm}b\eta(a)=[(\rm
id\otimes \varphi)W]\eta(a)\end{eqnarray*} which proves \eqref{rown}
for $a=\mathcal{Q}(f_\omega)$. But the set
\[\{a=(\omega\otimes \mbox{id})W:f_\omega\in L^2(G)\}\]
 is a core for $\eta$, hence
equation \rf{rown} is satisfied for any $a\in\mathfrak{N}$.
\end{proof}
\begin{uw} Let  $\pi^R,\pi^L\in\Rep(\C^*(\Gamma),L^2(G))$
be representations that send generators
$u_\gamma\in\M(\C^*(\Gamma))$ to $R_\gamma$ and
$L_\gamma\in\B(L^2(G))$ respectively. Let $f,\tilde
f\in\M(\C^*(\Gamma))$ and  $\omega\in\B(L^2(G))_*$ be such  that
$f_\omega\in L^2(G)$. Then using a method similar to the one used
in the proof of Proposition \ref{haar1} we can show that
\begin{equation}\label{fl} f_{\pi^R(f)\,\cdot\omega\cdot \pi^R(\tilde f)}=\pi^R(f)\pi^L(\kappa(\tilde f))(f_\omega)\end{equation}
where $\kappa$ is the coinverse on $\C^*(\Gamma)$.
\end{uw}
 With GNS-map $\eta:\mathfrak{N}\rightarrow L^2(G)$ we can associate a weight $h^\Psi:
\C_\infty(G)^{\widetilde{\Psi}\otimes\Psi}_+\mapsto \mathbb{R}_+$
: $h^\Psi(a^*a)=(\eta(a)|\eta(a))$.
\begin{stwr}\label{tr} Let $h^\Psi$ be the
 weight on $\C_\infty(G)^{\widetilde{\Psi}\otimes\Psi}$ introduced above.
Then it is a faithful trace. In particular it is strictly
faithful.\end{stwr}
\begin{proof} Let $u\in\M(\C^*(\Gamma))$ be the unitary element which appears in formula \eqref{uop}.
From the above remark and Proposition \ref{man} it follows that
\begin{eqnarray*}\eta(((\omega_{x,y}\otimes\id)W)^*)
&\hspace*{-0,25cm}=&\hspace*{-0,25cm}\eta((\pi^R(u)\cdot
\omega_{\bar x,\bar y}\cdot\pi^R(u)\otimes\id)W)\\
&\hspace*{-0,25cm}=&\hspace*{-0,25cm}\pi^R(u)\pi^L(\kappa(u))\eta((\omega_{\bar
x,\bar
y}\otimes\id)W)\\
&\hspace*{-0,25cm}=&\hspace*{-0,25cm}\pi^R(u)\pi^L(\kappa(u))\overline{\eta((\omega_{x,y}\otimes\id)W)}.\end{eqnarray*}
The set
\[\{a=(\omega\otimes \mbox{id})W:f_\omega\in L^2(G)\}\]
 is a core for $\eta$, hence we have
 \[\eta(a^*)=\pi^R(u)\pi^L(\kappa(u))\overline{\eta(a)}\]
 for any $a\in\mathfrak{N}$.
Now we can prove the trace property:
\begin{eqnarray*}h^\Psi(a^*a)&\hspace*{-0,25cm}=&\hspace*{-0,25cm}(\eta(a)|\eta(a))=(\overline{\eta(a)}|\overline{\eta(a)})
\\&\hspace*{-0,25cm}=&\hspace*{-0,25cm}(\pi^L(\kappa(u))^*\pi^R(u)^*\eta(a^*)|\pi^L(\kappa(u))^*\pi^R(u)^*\eta(a^*))\\
&\hspace*{-0,25cm}=&\hspace*{-0,25cm}(\eta(a^*)|\eta(a^*))=h^\Psi(aa^*).
\end{eqnarray*}
Let us prove the faithfulness of $h^\Psi$. Assume that
$h^\Psi(a^*a)=0$. Then \[h^\Psi(a^*c^*ca)=0=h^\Psi(caa^*c^*)\]
hence $\eta(a^*c^*)=a^*\eta(c^*)=0$. The set of elements of the
form $\eta(c^*)$ is dense in $L^2(G)$, hence $a=0$. The notion of
strict faithfulness was introduced in \cite{MNW}. It can be shown
that a faithful trace  is automatically strictly faithful. This
ends our proof.
\end{proof} Using Propositions \ref{haar1} and \ref{tr} one can
check that the assumptions of Theorem 3.9 of \cite{MNW} are
satisfied. Hence we get
\begin{twr} Let $( \C_\infty(G)^{\widetilde{\Psi}\otimes\Psi},\Delta^\Psi)$ be the
 quantum group with the multiplicative unitary $W$ and the weight $h^\Psi$ considered
 above. Then $h^\Psi$ is a Haar measure for $(
\C_\infty(G)^{\widetilde{\Psi}\otimes\Psi},\Delta^\Psi)$ and $W$
is the canonical multiplicative unitary.
\end{twr}
\end{subsection}
\end{section}
\begin{section}{An example of quantization of $SL(2,\mathbb{C})$.}
In this section we use the Rieffel deformation to quantize the
special linear group:
\[SL(2,\mathbb{C})=\left\{\begin{pmatrix}
  \alpha & \beta \\
  \gamma & \delta
\end{pmatrix}:\alpha,\beta,\gamma,\delta\in\mathbb{C},\,\alpha\delta-\beta\gamma=1\right\}.\]
(In what follows $SL(2,\mathbb{C})$ will be denoted by $G$.)
 The resulting quantum group is  the $\C^*$-algebraic version of
one of the  $*$-Hopf algebras introduced by  S.L. Woronowicz and
S. Zakrzewski in paper \cite{W9}. As a $*$-algebra it is generated
by four elements $\hat{\alpha},\hat{\beta},\hat\gamma,\hat\delta$,
satisfying the following commutation relations:
\[\begin{array}{rcl}
                    \hat{\alpha}\hat{\beta}&\hspace*{-0.25cm}=&\hspace*{-0.25cm}\hat{\beta}\hat{\alpha}\\
                    \hat{\alpha}\hat{\delta}&\hspace*{-0.25cm}=&\hspace*{-0.25cm}\hat{\delta}\hat{\alpha}\\
                    \hat{\alpha}\hat{\gamma}&\hspace*{-0.25cm}=&\hspace*{-0.25cm}\hat{\gamma}\hat{\alpha}\\
                    \hat{\beta}\hat{\gamma}&\hspace*{-0.25cm}=&\hspace*{-0.25cm}\hat{\gamma}\hat{\beta}\\
                    \hat{\beta}\hat{\delta}&\hspace*{-0.25cm}=&\hspace*{-0.25cm}\hat{\delta}\hat{\beta}\\
                    \hat{\gamma}\hat{\delta}&\hspace*{-0.25cm}=&\hspace*{-0.25cm}\hat{\delta}\hat{\gamma}\\
                    \hat{\alpha}\hat{\delta}&\hspace*{-0.25cm}=&\hspace*{-0.25cm} 1+\hat{\beta}\hat{\gamma}

\end{array}
\]\begin{equation}\label{relhl}
\begin{array}{rclrclrccrcl}\hat{\alpha}\hat{\alpha}^*&\hspace*{-0.25cm}=
                    &\hspace*{-0.25cm}\hat{\alpha}^*\hat{\alpha}&&&&&&&&&\\
                    \hat{\alpha}\hat{\beta}^*&\hspace*{-0.25cm}=
                    &\hspace*{-0.25cm}t\hat{\beta}^*\hat{\alpha}&\hat{\beta}\hat{\beta}^*&\hspace*{-0.25cm}=
                    &\hspace*{-0.25cm}\hat{\beta}^*\hat{\beta}&&&&&&\\
                    \hat{\alpha}\hat{\gamma}^*&\hspace*{-0.25cm}=
                    &\hspace*{-0.25cm}t^{-1}\hat{\gamma}^*\hat{\alpha}&\hat{\beta}\hat{\gamma}^*&\hspace*{-0.25cm}=
                    &\hspace*{-0.25cm}\hat{\gamma}^*\hat{\beta}&\hat{\gamma}\hat{\gamma}^*&\hspace*{-0.25cm}=
                    &\hspace*{-0.25cm}\hat{\gamma}^*\hat{\gamma}&&&\\
                    \hat{\alpha}\hat{\delta}^*&\hspace*{-0.25cm}=
                    &\hspace*{-0.25cm}\hat{\delta}^*\hat{\alpha}&\hat{\beta}\hat{\delta}^*&\hspace*{-0.25cm}=
                    &\hspace*{-0.25cm}t^{-1}\hat{\delta}^*\hat{\beta}&\hat{\gamma}\hat{\delta}^*&\hspace*{-0.25cm}=
                    &\hspace*{-0.25cm}t\hat{\delta}^*\hat{\gamma}&\hat{\delta}\hat{\delta}^*&\hspace*{-0.25cm}=
                    &\hspace*{-0.25cm}\hat{\delta}^*\hat{\delta}
\end{array} \end{equation} where $t$ is a nonzero real parameter.
The comultiplication, coinverse and counit act on them in the standard
way:
\begin{equation}\label{kodzialh}\begin{array}{ccc}
\begin{array}{rcl}
\Delta(\hat{\alpha})&\hspace*{-0.25cm}=&\hspace*{-0.25cm}\hat{\alpha}\otimes\hat{\alpha}+\hat{\beta}\otimes\hat{\gamma}\\
\Delta(\hat{\beta})&\hspace*{-0.25cm}=&\hspace*{-0.25cm}\hat{\alpha}\otimes\hat{\beta}+\hat{\beta}\otimes\hat{\delta}\\
\Delta(\hat{\gamma})&\hspace*{-0.25cm}=&\hspace*{-0.25cm}\hat{\gamma}\otimes\hat{\alpha}+\hat{\delta}\otimes\hat{\gamma}\\
\Delta(\hat{\delta})&\hspace*{-0.25cm}=&\hspace*{-0.25cm}\hat{\gamma}\otimes\hat{\beta}+\hat{\delta}\otimes\hat{\delta}
\end{array}&

\begin{array}{rcl}
\kappa(\hat{\alpha})&\hspace*{-0.25cm}=&\hspace*{-0.25cm}\hat{\delta}\\
\kappa(\hat{\beta})&\hspace*{-0.25cm}=&\hspace*{-0.25cm}-\hat{\beta}\\
\kappa(\hat{\gamma})&\hspace*{-0.25cm}=&\hspace*{-0.25cm}-\hat{\gamma}\\
\kappa(\hat{\delta})&\hspace*{-0.25cm}=&\hspace*{-0.25cm}\hat{\alpha}
\end{array}&
\begin{array}{rcl}
\varepsilon(\hat{\alpha})&\hspace*{-0.25cm}=&\hspace*{-0.25cm}1\\
\varepsilon(\hat{\beta})&\hspace*{-0.25cm}=&\hspace*{-0.25cm}0\\
\varepsilon(\hat{\gamma})&\hspace*{-0.25cm}=&\hspace*{-0.25cm}0\\
\varepsilon(\hat{\delta})&\hspace*{-0.25cm}=&\hspace*{-0.25cm}1.
\end{array}
\end{array}
\end{equation}

The deformation procedure in our example is based on the abelian
subgroup $\Gamma\subset G$ of diagonal matrices:
\begin{displaymath} \Gamma=\left\{\left(\begin{array}{cc}
  w & 0 \\
  0 & w^{-1}\end{array}\right):w\in
  \mathbb{C}_*\right\}.
\end{displaymath}
To simplify some calculations we pull back the action of $\Gamma^2$
on $\C_\infty(G)$ to the action  of $\mathbb{C}^2$ on
$\C_\infty(G)$. The resulting action is denoted by $\rho$:
\begin{equation}\label{conrho}
\bigl(\rho_{z_1,z_2}f\bigr)(g)=f\left(\left(\begin{array}{cc}
  e^{-z_1} & 0 \\
  0 & e^{z_1}\end{array}\right)g\left(\begin{array}{cc}
  e^{z_2} & 0 \\
  0 & e^{-z_2}\end{array}\right)\right)
\end{equation}
 Let us fix a $2$-cocycle on the dual group.
The additive group $(\mathbb{C},+)$ is self dual, with the duality
given by:
\[\mathbb{C}^2\ni(z_1,z_2)\mapsto\exp\bigl({\ir\rm
Im}(z_1z_2)\bigr)\in \mathbb{T}.\] Let $s\in \mathbb{R}$. For any
$z_1,z_2\in \mathbb{C}$ we set
\[\Psi(z_1,z_2)=\exp\bigl(\ir s{\rm Im}(z_1\bar{z}_2)\bigr).\]
It is clear, that $\Psi\in\C_b(\mathbb{C}^2)$ satisfies the
$2$-cocycle condition. Using results of Section \ref{qg} we deform
the standard $\mathbb{C}^2$-product structure on
$\C_\infty\bigl(G\bigr)\rtimes_\rho \mathbb{C}^2$
 to $\Bigl(\C_\infty\bigl(G\bigr)\rtimes_\rho \mathbb{C}^2,\lambda,\hat{\rho}^{\widetilde{\Psi}\otimes\Psi}\Bigl)$.
 In our case $\widetilde{\Psi}$ is just the complex conjugate of $\Psi$
and the deformed action of the dual group is given by
\[
\hat{\rho}^{\widetilde{\Psi}\otimes\Psi}_{z_1,z_2}(b)= \lambda_{-s\bar{z}_1,s\bar{z}_2}\hat{\rho}_{z_1,z_2}(b)
 \lambda_{-s\bar{z}_1,s\bar{z}_2}^* \] for any $b\in \C_\infty\bigl(G\bigr)\rtimes_\rho
\mathbb{C}^2$. The  Landstad algebra $A$ of the deformed
$\mathbb{C}^2$-product
 carries the structure of a quantum group. Our aim is to show that this
quantum group is the $\C^*$-algebraic version of the Hopf
$*$-algebra described above. The relation between parameters
$s,t\in\mathbb{R}$ is $t=\exp(-2s)$.
\begin{subsection}{$\C^*$-algebra structure.}
In this section we will construct four affiliated elements $\hat{\alpha},\hat{\beta},\hat{\gamma},\hat{\delta}\,\,\eta\,\, A$
 and show that they generate  $\C^*$-algebra $A$.

Let $T_r,T_l\in\C^*(
\mathbb{C}^2)^\eta\subset\Bigl(\C_\infty\bigl(G\bigr)\rtimes_\rho
\mathbb{C}^2\Bigr)^\eta$ be infinitesimal generators of the left
and right shifts. By definition $T_l$ and $T_r$ are normal
elements satisfying:
\begin{equation}\label{inf} \lambda_{z_1,z_2}=\exp\bigl(\ir{\rm Im}(z_1T_l)\bigr)\exp\bigl(\ir{\rm Im}(z_2T_r)\bigr)\end{equation}
for any $z_1,z_2\in\mathbb{C}$. Let $\alpha,\beta,\gamma,\delta$ be
coordinate functions on $G$:
 \[\alpha,\beta,\gamma,\delta\in\C\bigl(G\bigr)=
 \Bigl(\C_\infty\bigl(G\bigr)\Bigr)^\eta\subset\Bigl(\C_\infty\bigl(G\bigr)\rtimes_\rho
\mathbb{C}^2\Bigr)^\eta.\]
 Consider also a unitary element:
\begin{equation} \label{udef} U=\exp\bigl(\ir s{\rm Im}(T_r^*T_l)\bigr)\in\M\bigl(\C^*(
\mathbb{C}^2)\bigr)\subset \M\bigl(
\C_\infty\bigl(G\bigr)\rtimes_\rho \mathbb{C}^2\bigr).\end{equation}
We use it to define four normal elements affiliated with
$\C_\infty\bigl(G\bigr)\rtimes_\rho \mathbb{C}^2$:
\begin{equation}\label{formgen}\begin{array}{cc}
 \hat{ \alpha}=U\alpha U^* & \hat{ \beta}=U^*\beta U \\
 \hat{\gamma}=U^*\gamma U &  \hat{ \delta}=U\delta U^*.
\end{array}
\end{equation}
In the next lemma we present different formulas for
$\hat\alpha,\hat\beta,\hat\gamma,\hat\delta$ which will be needed
later.
\begin{lem}\label{express}
Let
$\alpha,\beta,\gamma,\delta\in\C_\infty\bigl(G\bigr)^\eta\subset\Bigl(\C_\infty\bigl(G\bigr)\rtimes_\rho
\mathbb{C}^2\Bigr)^\eta $ be coordinate functions on $G$. Let
$T_l,T_r$ be infinitesimal generators defined
 by \eqref{inf} and let
\[\hat{ \alpha},\hat{ \beta},\hat{\gamma},\hat{
\delta}\in\Bigl(\C_\infty\bigl(G\bigr)\rtimes_\rho
\mathbb{C}^2\Bigr)^{\eta}\] be normal elements \eqref{formgen}. Then
\begin{equation}\label{expnew}\left\{
\begin{array}{llll}
1.&\alpha \mbox{ and } T_l+T_r\mbox{ strongly commute and }
\hat{\alpha}&\hspace*{-.25cm}=&\hspace*{-.25cm}\exp(-s(T_l^*+T_r^*))\alpha;\\
2.&\beta \mbox{ and } T_l-T_r\mbox{ strongly commute and }
\hat{\beta}&\hspace*{-.25cm}=&\hspace*{-.25cm}\exp(s(T_l^*-T_r^*))\beta;\\
3.&\gamma \mbox{ and }T_l-T_r\mbox{ strongly commute and }
\hat{\gamma}&\hspace*{-.25cm}=&\hspace*{-.25cm}\exp(s(T_r^*-T_l^*))\gamma;\\
4.&\delta \mbox{ and } T_l+T_r\mbox{ strongly commute and }
\hat{\delta}&\hspace*{-.25cm}=&\hspace*{-.25cm}\exp(s(T_l^*+T_r^*))\delta.
\end{array}\right.
\end{equation}
\end{lem}
\begin{proof}
The fact that  $T_l+T_r$ and $\alpha$ strongly commute  follows
from the identity
\[\exp\bigl(\ir{\rm Im}(z(T_l+T_r))\bigr)\alpha\exp\bigl(-\ir{\rm
Im}(z(T_l+T_r))\bigr)=\alpha.\] We check it below:
\[\exp\bigl(\ir{\rm Im}(z(T_l+T_r))\bigr)\alpha\exp\bigl(-\ir{\rm
Im}(z(T_l+T_r))\bigr)=\lambda_{z,z}\alpha\lambda_{z,z}^*=\exp(-z+z)\alpha=\alpha.\]
To prove the equality $\hat{\alpha}=\exp(-s(T_l^*+T_r^*))\alpha$
note that
\begin{equation}\label{es}\begin{array}{rcl}\hat{\alpha}&\hspace*{-0,25cm}=&\hspace*{-0,25cm}\exp\bigl(\ir s{\rm
Im}(T_r^*T_l^{})\bigr)\alpha\exp\bigl(-\ir s{\rm
Im}(T_r^*T_l^{})\bigr)\\&\hspace*{-0,25cm}=&\hspace*{-0,25cm}\exp\bigl(\ir
s{\rm Im} ((T_l+T_r)^*T_l)\bigr)\alpha\exp\bigl(-\ir s{\rm
Im}((T_l+T_r)^*T_l)\bigr),\end{array}\end{equation} where  we used
the fact that $\exp\bigl(\ir s{\rm Im}(T_l^*T_l^{})\bigr)=1$.
Using the strong commutativity of
 $T_l+T_r$ and $\alpha$ and  the following identity:
\[\exp\bigl(\ir s{\rm Im}(wT_l)\bigr)\alpha\exp\bigl(-\ir s{\rm
Im}(wT_l)\bigr)=\exp(-sw)\alpha,\] we get
$\hat{\alpha}=\exp(-s(T_l^*+T_r^*))\alpha$. This ends the proof of
 point 1 of \eqref{expnew}.
 Using the same techniques we prove points 2,3,4.
\end{proof} Our objective is to show that
$\hat{\alpha},\hat{\beta}, \hat{\gamma}, \hat{\delta}$ are
generators of  $\C^*$-algebra $A$. In particular we have to show
that they are affiliated with $A$. The following proposition is
the first step toward the proof of this fact.
\begin{stwr}
\label{eg2} Let $\Bigl(\C_\infty\bigl(G\bigr)\rtimes_\rho
\mathbb{C}^2,\lambda,\hat{\rho}^{\widetilde{\Psi}\otimes\Psi}\Bigl)$
be the deformed $\mathbb{C}^2$- product, $A$  its Landstad algebra
and $\hat{\alpha}\in \left(\C_\infty\bigl(G\bigr)\rtimes_\rho
\mathbb{C}^2\right)^\eta$ the normal element defined in
\eqref{formgen}. Then $f(\hat{ \alpha})\in \M(A)$ for any $f\in
\C_\infty (\mathbb{C})$.
\end{stwr}
\begin{proof}
Let us first prove the invariance of  $f(\hat{ \alpha})$ under the
action $\hat{\rho}^{\widetilde{\Psi}\otimes\Psi}$. It is enough to
check that $\hat{ \alpha}$ is invariant. In order to do that we
calculate
\begin{equation}\label{intcal1}\lambda_{z_1,z_2}\alpha\lambda_{z_1,z_2}^*=
\rho_{z_1,z_2}(\alpha)=\exp(-z_1+z_2)\alpha.\end{equation}
Furthermore
\begin{eqnarray*} \hat{\rho}^{\widetilde{\Psi}\otimes\Psi}_{z_1,z_2}(\hat{
\alpha})&\hspace*{-.25cm}=&\hspace*{-.25cm}
\hat{\rho}^{\widetilde{\Psi}\otimes\Psi}_{z_1,z_2}(U\alpha
U^*)\\&\hspace*{-.25cm}=&\hspace*{-.25cm}
\hat{\rho}^{\widetilde{\Psi}\otimes\Psi}_{z_1,z_2}(U)\hat{\rho}^{\widetilde{\Psi}\otimes\Psi}_{z_1,z_2}
(\alpha)\hat{\rho}^{\widetilde{\Psi}\otimes\Psi}_{z_1,z_2}(U)^*.\end{eqnarray*}
We compute $\hat{\rho}^{\widetilde{\Psi}\otimes\Psi}_{z_1,z_2}(U)$
and $ \hat{\rho}^{\widetilde{\Psi}\otimes\Psi}_{z_1,z_2}(\alpha)$
separately:
\begin{eqnarray*}\hat{\rho}^{\widetilde{\Psi}\otimes\Psi}_{z_1,z_2}(\alpha)&\hspace*{-.25cm}=&\hspace*{-.25cm}
\lambda_{-s\bar{z}_1,s\bar{z}_2}\hat{\rho}_{z_1,z_2}(\alpha)
\lambda_{-s\bar{z}_1,s\bar{z}_2}^*\\&\hspace*{-.25cm}=&\hspace*{-.25cm}\lambda_{-s\bar{z}_1,s\bar{z}_2}\alpha
\lambda_{-s\bar{z}_1,s\bar{z}_2}^*\\&\hspace*{-.25cm}=&\hspace*{-.25cm}
\exp(s\bar{z}_1+s\bar{z}_2)\alpha\\\hat{\rho}^{\widetilde{\Psi}\otimes\Psi}_{z_1,z_2}(U)
&\hspace*{-.25cm}=&\hspace*{-.25cm}\hat{\rho}^{\widetilde{\Psi}\otimes\Psi}_{z_1,z_2}\bigl(\exp\bigl(\ir
s {\rm Im}(T_r^*T_l^{})\bigr)\bigr)
\\&\hspace*{-.25cm}=&\hspace*{-.25cm}\hat{\rho}_{z_1,z_2}\bigl(\exp\bigl(\ir
s{\rm Im}(T_r^*T_l^{})\bigr)\bigr)
\\&\hspace*{-.25cm}=&\hspace*{-.25cm}\exp\bigl(\ir
s{\rm Im}\bigl((T_r^*+\bar{z}_2)(T_l+z_1)\bigr)\bigr)
\\&\hspace*{-.25cm}=&\hspace*{-.25cm}U\lambda_{s\bar{z}_2,-s\bar{z}_1}\Psi(z_1,z_2).\end{eqnarray*}
Using \eqref{intcal1} we get
\begin{eqnarray*}\hat{\rho}^{\widetilde{\Psi}\otimes\Psi}_{z_1,z_2}(\hat{
\alpha})
&\hspace*{-.25cm}=&\hspace*{-.25cm}\exp(s\bar{z}_1+s\bar{z}_2)U\lambda_{s\bar{z}_2,-s\bar{z}_1}\alpha
\lambda_{-s\bar{z}_2,s\bar{z}_1}^*U^*
\\&\hspace*{-.25cm}=&\hspace*{-.25cm}\exp(s\bar{z}_1+s\bar{z}_2)\exp(-s\bar{z}_1-s\bar{z}_2)U\alpha
U^* =\hat{\alpha}.\end{eqnarray*} Let us now check that the map
\begin{equation}\label{2cont}\mathbb{C}^2\ni(z_1,z_2)\mapsto \lambda_{z_1,z_2}f(\hat{
\alpha})\lambda_{z_1,z_2}^*\in \M\bigl( \C_\infty(G)\rtimes
\mathbb{C}^2 \bigr)
\end{equation} is norm continuous. For this note that:
\begin{eqnarray*}\lambda_{z_1,z_2}f(\hat{
\alpha})\lambda_{z_1,z_2}^*&\hspace*{-.25cm}=&\hspace*{-.25cm}U\lambda_{z_1,z_2}f(
\alpha)\lambda_{z_1,z_2}^*U^*\\&\hspace*{-.25cm}=&\hspace*{-.25cm}Uf(e^{-z_1+z_2}\alpha)U^*.\end{eqnarray*}
Function $f$ is continuous and vanishes at infinity, hence we get
 norm continuity \eqref{2cont}. This shows that $f(\hat{\alpha})$
satisfies  the first and  second Landstad condition of \rf{lc1}
which is enough to be an element of $\M(A)$.
\end{proof}
To prove that $\hat{\alpha}$ is affiliated to $A$ we need one more
\begin{stwr}\label{eg3}
The set \[\mathcal{I}=\left\{f(\hat{\alpha})A: f\in
 \C_{\infty}( \mathbb{C})\right\}\]  is linearly dense in
$A$.
\end{stwr}
\begin{proof} Recall that $\rho^{\widetilde{\Psi}\otimes\Psi} $ is the action of $\mathbb{C}^2$ on $A$ implemented by
unitary elements $\lambda_{z_1,z_2}$.
 It is easy to see that $\mathcal{I}$ is invariant under
 $\rho^{\widetilde{\Psi}\otimes\Psi}$. Let $g\in\C_\infty(\mathbb{C})$ be a function given by the formula
 $g(z)=(1+\bar{z}z)^{-1}$.
Then $g(\hat\alpha)=U(1+\alpha^*\alpha)U^*$ and we have:
\begin{equation}\label{zawid}
\begin{array}{rcl}\left(\C^*( \mathbb{C}^2)g(\hat\alpha)A\C^*(
\mathbb{C}^2)\right)^{\rm
cls}&\hspace{-0,25cm}=&\hspace{-0,25cm}\left(\C^*(
\mathbb{C}^2)(1+\alpha^*\alpha)^{-1}U^*A\C^*(
\mathbb{C}^2)\right)^{\rm
cls}\\&\hspace{-0,25cm}\subset&\hspace{-0,25cm}\left(\C^*(
\mathbb{C}^2)\mathcal{I}\C^*( \mathbb{C}^2)\right)^{\rm
cls}\end{array} \end{equation} where we used the equality $\C^*(
\mathbb{C}^2)U=\C^*( \mathbb{C}^2)$. Note that the set $U^*A\C^*(
\mathbb{C}^2)$ is linearly dense in $
\C_\infty\bigl(G\bigr)\rtimes_\rho \mathbb{C}^2$. Using the fact
that $\alpha$ is affiliated with
 $ \C_\infty\bigl(G\bigr)\rtimes_\rho \mathbb{C}^2$ we see that the set
 $\C^*( \mathbb{C}^2)(1+\alpha^*\alpha)^{-1}U^*A\C^*( \mathbb{C}^2)$ is linearly dense in
  $ \C_\infty\bigl(G\bigr)\rtimes_\rho \mathbb{C}^2$. Hence by
\eqref{zawid} the set  $\C^*( \mathbb{C}^2)\mathcal{I}\C^*(
\mathbb{C}^2)$ is linearly dense in
$\C_\infty\bigl(G\bigr)\rtimes_\rho \mathbb{C}^2$. Using Lemma
\ref{dsw} we get the linear density of $\mathcal{I}$ in $A$.
\end{proof}
Let us define the homomorphism of $\C^*$-algebras:
\[  \C_\infty (\mathbb{C})\ni f\mapsto \pi(f)= f(
\hat{\alpha})\in \M(A).\]
\begin{twr} Let $\pi$ be the homomorphism defined above. $\pi$
 is a morphism of $\C^*$-algebras: $\pi\in\Mor\bigl(\C_\infty (\mathbb{C});A\bigr)$.
 In particular $\hat{\alpha}$ is the normal
element affiliated with $A$.
\end{twr}
\begin{proof} By Proposition
\ref{eg3} we have
$\overline{\pi(\C_\infty(\mathbb{C}))A}^{\|\cdot\|}=A$ which shows
 that $\pi\in\Mor(\C_\infty(\mathbb{C});A)$. Let ${\rm
id}\in\C_\infty(\mathbb{C})^{\eta}$ be the identity function: ${\rm
id}(z)=z$ for all $z\in\mathbb{C}$. Applying morphism $\pi$ to ${\rm
id}\in\C_\infty(\mathbb{C})^{\eta}$ we get $\pi({\rm id})={\rm
id}(\hat\alpha)=\hat\alpha\in A^\eta$.
\end{proof} Using the same techniques we show that
$\hat{\beta}, \hat{\gamma}, \hat{\delta}\af A$. In the next
theorem we prove that they are in fact generators of $A$.
\begin{twr}\label{eg4} Let $\hat\alpha,\hat\beta,\hat\gamma,\hat\delta\aff A$
be affiliated elements introduced in \eqref{formgen}. Let us
consider the set:
\[\mathcal{V}=\left\{ f_1(
\hat{\alpha})f_2(\hat{\beta})f_3(\hat{\gamma})f_4( \hat{\delta}):
f_1,f_2,f_3,f_4\in \C_\infty (\mathbb{C})\right\}\subset\M(A).\]
Then $\mathcal{V}$ is a subset of $A$ and $\mathcal{V}^{\,\rm
cls}=A$. In particular $A$ is generated by elements $
\hat{\alpha},\hat{\beta},\hat{\gamma},\hat{\delta}\in A^\eta$.
\end{twr}
\begin{proof}
Let us start with a proof that $\mathcal{V}\subset A$. Mimicking
the proof of Theorem \ref{eg2} we show that elements of
$\mathcal{V}$ satisfy the first and the second Landstad condition
\eqref{lc1}. To check that they also satisfy the third one, we
need to show that
\begin{equation}\label{propd}x f_1(
\hat{\alpha})f_2(\hat{\beta})f_3(\hat{\gamma})f_4(
\hat{\delta})y\in \C_\infty\bigl(G\bigr)\rtimes_\rho
\mathbb{C}^2\end{equation}  for any $x,y\in \C^*( \mathbb{C}^2)$.
Let us consider the set
\[ \mathcal{W}=\{xf_1(
\hat{\alpha})f_2(\hat{\beta})f_3(\hat{\gamma})f_4(
\hat{\delta})y:f_1,f_2,f_3,f_4\in \C_\infty( \mathbb{C}),\,x,y\in
\C^*( \mathbb{C}^2)\}^{\rm cls}.\] Note that
$\mathcal{W}=(\C^*(\mathbb{C}^2)\mathcal{V}\C^*(\mathbb{C}^2))^{\rm
cls}$. We will show that: \[
\mathcal{W}=\C_\infty\bigl(G\bigr)\rtimes_\rho \mathbb{C}^2\] which
is a stronger property than \eqref{propd}.
 Using \eqref{formgen} we get
\begin{eqnarray*}\mathcal{W}=\bigl\{xUf_1(
\alpha){U^*}^2f_2(\beta)f_3(\gamma)U^2f_4(
\delta)U^*y:&\\&\hspace*{-2cm}f_1,f_2,f_3,f_4\in \C_\infty(
\mathbb{C}),\,x,y\in \C^*( \mathbb{C}^2)\bigr\}^{\rm
cls}.\end{eqnarray*} By unitarity of $U$ we can substitute $x$
with  $xU^*$ and  $y$ with $Uy$ not changing $\mathcal{W}$:
\[\mathcal{W}=\bigl\{xf_1( \alpha){U^*}^2f_2(\beta)f_3(\gamma)U^2f_4(
\delta)y:f_1,f_2,f_3,f_4\in \C_\infty( \mathbb{C}),\,x,y\in \C^*(
\mathbb{C}^2)\bigr\}^{\rm cls}.\] The map
\[\mathbb{C}^2\ni(z_1,z_2)\mapsto \rho_{z_1,z_2}(f(\alpha))=f(\exp(-z_1+z_2)\alpha)\] is norm
continuous, hence: \[\{f(\alpha)x:f\in\C_\infty(
\mathbb{C}),x\in\C^*( \mathbb{C}^2)\}^{\rm
cls}=\{xf(\alpha):f\in\C_\infty( \mathbb{C}),x\in\C^*(
\mathbb{C}^2)\}^{\rm cls}.\] In particular
\[\mathcal{W}= \bigl\{f_1(
\alpha)x{U^*}^2f_2(\beta)f_3(\gamma)U^2f_4(
\delta)y:f_1,f_2,f_3,f_4\in \C_\infty( \mathbb{C}),\,x,y\in \C^*(
\mathbb{C}^2)\bigr\}^{\rm cls}.\] Similarly, we commute $f_4(
\delta)$ and $y$: \[\mathcal{W}=\bigl\{f_1(
\alpha)x{U^*}^2f_2(\beta)f_3(\gamma)U^2yf_4(
\delta):f_1,f_2,f_3,f_4\in \C_\infty( \mathbb{C}),\,x,y\in \C^*(
\mathbb{C}^2)\bigr\}^{\rm cls}.\] Substituting $x$ with $xU^2$ and
$y$ with $U^{*2}y$ we get
\[\mathcal{W}= \{f_1( \alpha)xf_2(\beta)f_3(\gamma)yf_4(
\delta):f_1,f_2,f_3,f_4\in \C_\infty( \mathbb{C}),\,x,y\in \C^*(
\mathbb{C}^2)\}^{\rm cls}.\] Commuting back  $f_1( \alpha)$ ($f_4(
\delta)$ resp.) and $x$ ($y$ resp.) we obtain
\[\mathcal{W}=
\{xf_1( \alpha)f_2(\beta)f_3(\gamma)f_4(
\delta)y:f_1,f_2,f_3,f_4\in \C_\infty( \mathbb{C}),\,x,y\in \C^*(
\mathbb{C}^2)\}^{\rm cls}.\]
 The last set is obviously the whole $ \C_\infty\bigl(G\bigr)\rtimes_\rho \mathbb{C}^2$.
 Therefore we conclude that elements of $\mathcal{V}$ satisfies the Landstad conditions
 and  $\mathcal{V}\subset A$.
Moreover $\mathcal{V}$ is
$\rho^{\widetilde{\Psi}\otimes\Psi}$-invariant and the set $\C^*(
\mathbb{C}^2)\mathcal{V}\C^*( \mathbb{C}^2)$ is linearly dense in
$\C_\infty\bigl(G\bigr)\rtimes_\rho \mathbb{C}^2$. Using Lemma
\ref{dsw} we see that $\mathcal{V}^{\rm cls}=A$. In particular
$\hat{\alpha},\hat{\beta},\hat{\gamma},\hat{\delta}$ separate
representations of $A$ and
\[(1+\hat{\alpha}^*\hat{\alpha})^{-1}(1+\hat{\beta}^*\hat{\beta})^{-1}(1+\hat{\gamma}^*\hat{\gamma})^{-1}(1+\hat{\delta}^*\hat{\delta})^{-1}\in
A.\] By Theorem 3.3 of \cite{W4} we see that $A$ is generated by
$\hat\alpha,\hat\beta,\hat\gamma,\hat\delta$.
\end{proof}
\end{subsection}
\begin{subsection}{Commutation relations.}
The aim of this section is to show that generators
$\hat{\alpha},\hat{\beta},\hat{\gamma},\hat{\delta}$ satisfy
relations \eqref{relhl}. Note that in general it is impossible to
multiply affiliated elements, so we  have to give a precise
meaning to \eqref{relhl}. We start with  considering a more
general type of relations.
 Let $p,q$ be real, strictly positive numbers and $(R,S)$ a pair of normal operators
 acting on  $H$.
 The precise meaning of the relations
\begin{eqnarray*}
RS&\hspace*{-.25cm}=&\hspace*{-.25cm}pSR\\RS^*&\hspace*{-.25cm}=&\hspace*{-.25cm}qS^*R.\end{eqnarray*} was given in \cite{W3}:
\begin{defin} \label{mu1}Let $(R,S)$ be a pair of normal operators acting on a Hilbert
 space $H$. We say that $(R,S)$ is a $(p,q)$-commuting pair if
\begin{itemize}
\item[1.] $|R|$ and $|S|$ strongly commute.
\item[2.] $(\mbox{Phase}R)(\mbox{Phase}S)=(\mbox{Phase}S)(\mbox{Phase}R)$.
\item[3.] On $\ker R^{\perp}$ we have
\[(\mbox{Phase}R)|S|(\mbox{Phase}R)^*=\sqrt{pq}\,|S|.\]
\item[4.]  On $\ker S^{\perp}$ we have
\[(\mbox{Phase}S)|R|(\mbox{Phase}S)^*=\sqrt{q/p}\,|R|.\]
\end{itemize}
The set of all $(p,q)$-commuting pairs of normal operators acting on
a Hilbert space $H$ is denoted by $D_{p,q}(H)$. Note that
$(1,1)$-commuting pair of normal operators is just a strongly
commuting pair of operators.
\end{defin}

We need a version of the above definition which is suitable for a
pair of normal elements affiliated with a $\C^*$-algebra. In what
follows we shall use  the symbol $z(T)$  to denote the
$z$-transform of an element $T$: $z(T)=T(1+T^*T)^{-\frac{1}{2}}$.
\begin{defin} \label{mu2} Let $A$ be a $\C^*$-algebra and $(R,S)$
a pair of normal elements affiliated with $A$. We say that $(R,S)$
is a $(p,q)$-commuting pair if
\begin{itemize}
\item[1.]$z(R)z(S^*)=z(\sqrt{pq}\,S^*)z(\sqrt{q/p}\,R)$
\item[2.]$z(\sqrt{q/p}\, R)z(S)=z(\sqrt{pq}\, S)z(R)$.
\end{itemize}
The set of all $(p,q)$-commuting pairs of normal elements affiliated with a $\C^*$-algebra $A$ is denoted by $D_{p,q}(A)$.
\end{defin}
It turns out that  Definitions \ref{mu1} and \ref{mu2}  are in a
sense equivalent. Namely we have:
\begin{stwr}\label{eqiu} Let $(R,S)$ be a pair of normal operators acting on $H$. It is a $(p,q)$-commuting pair in the sense of Definition
\ref{mu1} if and only if
\begin{equation}\label{pqcom}\left\{\begin{array}{rcl}
z(R)z(S^*)&\hspace*{-.25cm}=&\hspace*{-.25cm}z(\sqrt{pq}\,S^*)z(\sqrt{q/p}\,R)\\
z(\sqrt{q/p}\,R)z(S)&\hspace*{-.25cm}=&\hspace*{-.25cm}z(\sqrt{pq}\,S)z(R).
\end{array}\right.
\end{equation}
\end{stwr}
\begin{proof}
It is easy to see that a pair $(R,S)$ of $(p,q)$-commuting
operators satisfies \eqref{pqcom}. We will  prove the opposite
implication. Using  \eqref{pqcom} we get:
\begin{equation}\label{eqmd}\begin{array}{rcl}z({\sqrt{q/p}\, R})z(S)z(S)^*&\hspace*{-0.25cm}=&\hspace*{-0.25cm}z(\sqrt{pq}\, S)z(R)z(S)^*\\
&\hspace*{-0.25cm}=&\hspace*{-0.25cm}z(\sqrt{pq}\, S)z(\sqrt{pq}\,
S)^*z({\sqrt{q/p}\,R}).
\end{array}
\end{equation} Hence
\begin{equation}\label{eqmd1}\begin{array}{rcl}
z(\sqrt{q/p}\,R)^*z({\sqrt{q/p}\,
R})z(S)z(S)^*&&\\&\hspace*{-4cm}=&\hspace*{-2cm}
z(\sqrt{q/p}\,R)^*z(\sqrt{pq}\, S)z(\sqrt{pq}\,
S)^*z(\sqrt{q/p}\,R)\\&\hspace*{-4cm}=&\hspace*{-2cm}z(S)z(S)^*
z({\sqrt{q/p}\, R})^*z({\sqrt{q/p}\, R}).\end{array}\end{equation}
$T$ is a normal operator, hence $z_T^*z_T=z_Tz_T^*=z_{|T|}^2$  and
we get
\[\left(z(\sqrt{q/p}\,| R|)\right)^2z(|S|)^2=z(|S|)^2\left(z(\sqrt{q/p}\,| R|)\right)^2.\]
This shows that $(|R| ,|S|)$ is a pair of strongly commuting
operators.

 Using the polar decomposition of normal operators $R$ and $S$ we rewrite the second equation of \eqref{pqcom}:
\[\Ph(R)z(\sqrt{q/p}\,|R|)z(|S|)\Ph(S)=\Ph(S)z(\sqrt{pq}\,|S|)z(|R|)\Ph(R)
.\] Strong commutativity of $|R|$ and $|S|$ and  identities
\begin{eqnarray*}\Ph(S)\Ph(S)^*z(|S|)&\hspace*{-0.25cm}=&\hspace*{-0.25cm}z(|S|)\\
\Ph(R)\Ph(R)^*z(|R|)&\hspace*{-0.25cm}=&\hspace*{-0.25cm}z(|R|)
\end{eqnarray*} gives
\begin{eqnarray*}\Ph(R)\Ph(S)\Ph(S)^*z(|S|)z(\sqrt{q/p}\,|R|)\Ph(S)&\\&\hspace*{-6cm}=
\Ph(S)\Ph(R)\Ph(R)^*z(|R|)z(\sqrt{pq}\,|S|)\Ph(R).\end{eqnarray*}
Uniqueness of the polar decomposition implies that phases of $R$ and
$S$ commute:
\[\Ph(R)\Ph(S)=\Ph(S)\Ph(R).\]
Using equation \eqref{eqmd} we get
\[\Ph(R)z(\sqrt{q/p}\,|R|)z(|S|)=z(\sqrt{pq}\,|S|)\Ph(R)z(\sqrt{q/p}\,|R|).\] We already know
that $|R|$ and $|S|$ strongly commute, hence
\[\Ph(R)z(|S|)z(\sqrt{q/p}\,|R|)=z(\sqrt{pq}\,|S|)\Ph(R)z(\sqrt{q/p}\,|R|).\] This shows that on
$\ker R^{\perp}$ we have
\[\Ph(R)|S|\Ph(R)^*=\sqrt{pq}\,|S|.\] Similarly, one can prove that
$\Ph(S)|R|\Ph(S)^*=\sqrt{q/p}\,|R|$ on  $\ker S^{\perp}$.
\end{proof}
The next theorem shows that
$\hat\alpha,\hat\beta,\hat\gamma,\hat\delta\aff A$ satisfy relations
\eqref{relhl} in the sense of Definition \ref{mu2}.
\begin{twr}\label{comrel} Let  $
 \hat{ \alpha},\hat{ \beta},
 \hat{\gamma}, \hat{ \delta}\af A $
be elements given by \eqref{formgen}. Then
\begin{itemize}
\item[1.]
$\bigl(\hat{\alpha},\hat{\delta}\bigr),\bigl(\hat{\beta},\hat{\gamma}\bigr)\in
D_{1,1}(A)$
\item[2.] $\bigl(\hat{\alpha},\hat{\beta}\bigr)$,  $\bigl(\hat{\gamma},\hat{\delta}\bigr)\in D_{1,t}(A)$
\item[3.] $\bigl(\hat{\alpha},\hat{\gamma}\bigr)$, $\bigl(\hat{\beta},\hat{\delta}\bigr)\in D_{1,t^{-1}}(A)$
\end{itemize} where $t=\exp(-2s)$.
Consider normal elements $\hat{\alpha}\hat{\delta}\,,
\hat{\beta}\hat{\gamma}\af A$ (a product of two strongly commuting
normal elements is well defined). Then
$\bigl(\hat{\alpha}\hat{\delta},\hat{\beta}\hat{\gamma}\bigr)\in
D_{1,1}(A)$ and
\[\hat{\alpha}\hat{\delta}-\hat{\beta}\hat{\gamma}=1.\]
\end{twr}
\begin{proof} Directly from \eqref{formgen} it follows that $\bigl(\hat{\alpha},\hat{\delta}\bigr)\in
D_{1,1}\bigl(A\bigr)$ and $\bigl(\hat{\beta},\hat{\gamma}\bigr)\in
D_{1,1}\bigl(A\bigr)$. Note that the affiliated element
$\alpha\delta\,\eta \C_\infty(G)$ is $\rho$-invariant:
$\rho_{z_1,z_2}(\alpha\delta)=\alpha\delta$ where $\rho$ is the
action defined by \eqref{conrho}. Therefore, at the level of the
crossed product, $\alpha\delta$ commutes with
$\C^*(\mathbb{C}^2)$. Using the fact that
$U\in\M(\C^*(\mathbb{C}^2))$ we get
\[\hat{\alpha}\hat{\delta}=U\alpha\delta U^*=\alpha\delta.\] Similar
reasoning shows that $\hat{\beta}\hat{\gamma}=\beta\gamma$.
Therefore
$\bigl(\hat{\alpha}\hat{\delta},\hat{\beta}\hat{\gamma}\bigr)\in
D_{1,1}(A)$ and
\[\hat{\alpha}\hat{\delta}-\hat{\beta}\hat{\gamma}=\alpha\delta-\beta\gamma=1.\]
Now let us prove that $\bigl(\hat{\alpha},\hat{\beta}\bigr)\in
D_{1,t}\bigl(A\bigr)$. Using the faithful representation $\pi^{\rm
can}$ of $A$ on $L^2(G)$ we can treat  generators
$\hat{\alpha},\hat{\beta}$ as normal operators acting on $L^2(G)$.
We will show that
 $\bigl(\hat{\alpha},\hat{\beta}\bigr)\in
D_{1,t}\bigl(L^2(G)\bigr)$ which by Proposition \ref{eqiu} is
equivalent with the containment
$\bigl(\hat{\alpha},\hat{\beta}\bigr)\in D_{1,t}\bigl(A\bigr)$.
Using Lemma \ref{express} we get:
\begin{equation}\label{biega}\left\{\begin{array}{rcl}\mbox{Phase}(\hat{
\alpha})&\hspace*{-.25cm}=&\hspace*{-.25cm}\exp\bigl(\ir s{\rm
Im}(T_l+T_r)\bigr)\mbox{Phase}(\alpha)\\|\hat{\alpha}|&\hspace*{-.25cm}=&\hspace*{-.25cm}\exp\bigl(-s{\rm
Re}(T_l+T_r)\bigr)|\alpha|,\\
\mbox{Phase}(\hat{\beta})&\hspace*{-.25cm}=&\hspace*{-.25cm}
\exp\bigl(-\ir s{\rm
Im}(T_l-T_r)\bigr)\mbox{Phase}(\beta)\\|\hat{\beta}|&\hspace*{-.25cm}=&\hspace*{-.25cm}\exp\bigl(s{\rm
Re}(T_l-T_r)\bigr)|\beta|.\end{array}\right.\end{equation}
Moreover, it is easy to check that
\begin{equation}\label{biegc}\left\{\begin{array}{rcl}\exp\bigl(\ir s{\rm
Im}(T_l+T_r)\bigr)\beta\exp\bigl(-\ir s{\rm
Im}(T_l+T_r)\bigr)&\hspace*{-.25cm}=&\hspace*{-.25cm}\exp(-2s)\beta\\\exp\bigl(\ir
s {\rm Im}(T_l-T_r)\bigr)\alpha\exp\bigl(-\ir s{\rm
Im}(T_l-T_r)\bigr)&\hspace*{-.25cm}=&\hspace*{-.25cm}\exp(-2s)\alpha\\\exp\bigl(-\ir
s {\rm Re}(T_l+T_r)\bigr)|\beta|\exp\bigl(\ir s{\rm
Re}(T_l+T_r)\bigr)&\hspace*{-.25cm}=&\hspace*{-.25cm}|\exp(2\ir s
)\beta|=|\beta|\\\exp\bigl(-\ir s{\rm
Re}(T_l-T_r)\bigr)|\alpha|\exp\bigl(\ir s{\rm
Re}(T_l-T_r)\bigr)&\hspace*{-.25cm}=&\hspace*{-.25cm}|\exp(2\ir
s)\alpha|=|\alpha|.\end{array}\right.\end{equation} Equations
\eqref{biega} and \eqref{biegc}  show together that:
 \[ \begin{array}{ll}
1.&\mbox{Phase}(\hat{\alpha})\mbox{Phase}(\hat{\beta})=\mbox{Phase}(\hat{\beta})\mbox{Phase}(\hat{\alpha})\\
2.&\mbox{Phase}(\hat{\alpha})|\hat{\beta}|\mbox{Phase}(\hat{\alpha})^*=\exp(-2s)|\hat{\beta}|\\
3.&\mbox{Phase}(\hat{\beta})|\hat{\alpha}|\mbox{Phase}(\hat{\beta})^*=\exp(-2s)|\hat{\alpha}|\\
4.&|\hat{\alpha}|\mbox{ and } |\hat{\beta}| \mbox{ strongly
commute.}
 \end{array}\] Note that   $\ker\hat{\alpha}=\ker\hat{\beta}=\{0\}$ hence
$\bigl(\hat{\alpha},\hat{\beta}\bigr)\in
D_{1,t}\bigl(L^2(G)\bigr)$. Using the same techniques we prove all
other assertions of our theorem.
\end{proof}
\end{subsection}
\begin{subsection}{Comultiplication.} Let $\Delta^\Psi\in \Mor(A;A\otimes
A)$ be the  comultiplication on $A$. As was shown in Theorem
\ref{forcom}, it is given by:
\begin{equation}\label{comex}\Delta^\Psi(a)=\Upsilon\Delta(a)\Upsilon^*,\end{equation}
where $\Delta\in \Mor\bigl( \C_\infty(G)\rtimes \mathbb{C}^2;
\C_\infty(G)\rtimes \mathbb{C}^2\otimes
 \C_\infty(G)\rtimes \mathbb{C}^2\bigr)$
is uniquely characterized by two properties:
\begin{itemize}
\item
$\Delta(T_l)=T_l\otimes I, \,\,\Delta(T_r)=I\otimes T_r$;
 \item $\Delta$ restricted  to $ \C_\infty(G)$ coincides with
the comultiplication on $\C_\infty(G)$.
\end{itemize}
In our case the unitary element $\Upsilon$  is of the following
form: \begin{equation}\label{ups}\Upsilon=\exp\bigl(\ir s{\rm
Im}(T_r^*\otimes T_l)\bigr).\end{equation}
\begin{twr} Let $(A,\Delta^\Psi)$ be the quantum group considered above and
let $\hat{ \alpha}, \hat{ \beta}, \hat{\gamma},\hat{ \delta}$ be
the generators of $A$ given by \eqref{formgen}. Comultiplication
$\Delta^\Psi$ acts on generators in the standard way:
\begin{equation}
 \label{egcom}
\begin{array}{cc}
  \Delta^\Psi(\hat{\alpha})= \hat{\alpha}\otimes \hat{\alpha}+ \hat{\beta}\otimes \hat{\gamma}& \Delta^\Psi(\hat{\beta})=
   \hat{\alpha}\otimes \hat{\beta}+\hat{\beta}\otimes \hat{\delta} \\
   \Delta^\Psi(\hat{\gamma})= \hat{\gamma}\otimes \hat{\alpha}+ \hat{\delta}\otimes \hat{\gamma}& \Delta^\Psi(\hat{\delta})=
   \hat{\delta}\otimes \hat{\delta}+\hat{\gamma}\otimes \hat{\beta}.\end{array}
\end{equation}
\end{twr}
\begin{uw} The action of $\Delta^\Psi$ in the formula above is given
by the sum of  affiliated elements. In general it is not a well
defined operation. But in our case (as will be shown) this is a
sum of two normal strongly commuting  elements of $(A\otimes
A)^\eta$. This operation is well defined and gives a normal
element affiliated with $(A\otimes A)^\eta$.
\end{uw}
\begin{proof} Applying morphism $\Delta$ to $U\in\M(\C^*(\Gamma^2))$
(see  \eqref{udef}) we get:
\begin{equation}\label{ucom}\Delta(U)=\exp\bigl(\ir s{\rm Im}(T_l\otimes
T_r^*)\bigr).\end{equation} Let $R$ be a unitary element given by
the formula:
\[R=\exp\bigl(\ir s{\rm Im}(T_r^*\otimes T_l)\bigr)\exp\bigl(\ir s{\rm
Im}(T_l\otimes T_r^*)\bigr).\]  Using \eqref{comex}, \eqref{ups} and
\eqref{ucom} we get:
\begin{equation}\label{com2}
\begin{array}{lll}\Delta^\Psi(\hat{\alpha})&\hspace*{-.25cm}=&\hspace*{-.25cm}R(\alpha\otimes\alpha+\beta\otimes\gamma)R^*\\
&\hspace*{-.25cm}=&\hspace*{-.25cm}R(\alpha\otimes\alpha)R^*+
R(\beta\otimes\gamma)R^*.\end{array} \end{equation} Note that
 \begin{equation}\label{ur}(U^*\otimes U^*) R=\exp\bigl(\ir s{\rm Im}(T_r^*\otimes I-I\otimes T_r^*)(I\otimes
T_l-T_l\otimes I)\bigr).\end{equation} It is easy to check that
elements $(T_r\otimes I-I\otimes T_r)$ and $(I\otimes
T_l-T_l\otimes I)$ strongly commute with $\alpha\otimes \alpha$.
Hence by identity \eqref{ur}, $(U^*\otimes U^*) R$ commutes with
$\alpha\otimes \alpha$. Similarly, we check that the unitary
element $(U\otimes U) R$ commutes with $\beta\otimes \gamma$.
Using these two facts we get
\begin{equation}\label{com3}
\begin{array}{lll}R(\alpha\otimes\alpha)R^*&\hspace*{-.25cm}=&\hspace*{-.25cm}(U\otimes U)(U^*\otimes U^*)R(\alpha\otimes\alpha)R^*(U\otimes U)(U^*\otimes
U^*)\\&\hspace*{-.25cm}=&\hspace*{-.25cm}(U\otimes
U)(\alpha\otimes\alpha)(U^*\otimes
U^*)=\hat{\alpha}\otimes\hat{\alpha}\end{array}\end{equation} and
\begin{equation}\label{com4}
\begin{array}{lll}R(\beta\otimes\gamma)R^*&\hspace*{-.25cm}=&\hspace*{-.25cm}(U^*\otimes U^*)(U\otimes U)R(\beta\otimes\gamma)R^*(U^*\otimes
U^*)(U\otimes U)\\&\hspace*{-.25cm}=&\hspace*{-.25cm}(U^*\otimes
U^*)(\beta\otimes\gamma)(U\otimes U)=\hat{\beta}\otimes\hat{\gamma}.
\end{array}\end{equation} Equations \eqref{com2}, \eqref{com3},
\eqref{com4}  give:
\[\Delta^\Psi(\hat{\alpha})=
\hat{\alpha}\otimes\hat{\alpha}+\hat{\beta}\otimes\hat{\gamma}.\]
 All other assertions of our theorem are proven using the same
techniques.\end{proof}
\end{subsection}
\end{section}

\end{document}